\documentclass{amsart}
\begin{document}
\title[Overconvergence]{Differential overconvergence}
\author{Alexandru Buium and Arnab Saha}
\def \Rp{R_p}
\def \Rpi{R_{\pi}}
\def \dpi{\d_{\pi}}
\def \bT{{\bf T}}
\def \cI{{\mathcal I}}
\def \cJ{{\mathcal J}}
\def \ZN{\bZ[1/N,\zeta_N]}
\def \tA{\tilde{A}}
\def \o{\omega}
\def \tB{\tilde{B}}
\def \tC{\tilde{C}}
\def \alph{A}
\def \bet{B}
\def \bsigma{\bar{\sigma}}
\def \y{^{\infty}}
\def \Ra{\Rightarrow}
\def \uBS{\overline{BS}}
\def \lBS{\underline{BS}}
\def \lB{\underline{B}}
\def \<{\langle}
\def \>{\rangle}
\def \hL{\hat{L}}
\def \cU{\mathcal U}
\def \cF{\mathcal F}
\def \S{\Sigma}
\def \st{\stackrel}
\def \sd{Spec_{\d}\ }
\def \pd{Proj_{\d}\ }
\def \s{\sigma_2}
\def \i{\sigma_1}
\def \bs{\bigskip}
\def \cD{\mathcal D}
\def \cC{\mathcal C}
\def \cT{\mathcal T}
\def \cK{\mathcal K}
\def \cX{\mathcal X}
\def \sX{X_{set}}
\def \cY{\mathcal Y}
\def \cS{X}
\def \cR{\mathcal R}
\def \cE{\mathcal E}
\def \tcE{\tilde{\mathcal E}}
\def \cP{\mathcal P}
\def \cA{\mathcal A}
\def \cV{\mathcal V}
\def \cM{\mathcal M}
\def \cL{\mathcal L}
\def \cN{\mathcal N}
\def \tcM{\tilde{\mathcal M}}
\def \caS{\mathcal S}
\def \cG{\mathcal G}
\def \cB{\mathcal B}
\def \tG{\tilde{G}}
\def \cF{\mathcal F}
\def \h{\hat{\ }}
\def \hp{\hat{\ }}
\def \tS{\tilde{S}}
\def \tP{\tilde{P}}
\def \tA{\tilde{A}}
\def \tX{\tilde{X}}
\def \tcS{\tilde{X}}
\def \tT{\tilde{T}}
\def \tE{\tilde{E}}
\def \tV{\tilde{V}}
\def \tC{\tilde{C}}
\def \tI{\tilde{I}}
\def \tU{\tilde{U}}
\def \tG{\tilde{G}}
\def \tu{\tilde{u}}
\def \chu{\check{u}}
\def \tx{\tilde{x}}
\def \tL{\tilde{L}}
\def \tY{\tilde{Y}}
\def \d{\delta}
\def \e{\chi}
\def \bW{\mathbb W}
\def \bV{{\mathbb V}}
\def \bF{{\bf F}}
\def \bE{{\bf E}}
\def \bC{{\bf C}}
\def \bO{{\bf O}}
\def \bR{{\bf R}}
\def \bA{{\bf A}}
\def \bB{{\bf B}}
\def \cO{\mathcal O}
\def \ra{\rightarrow}
\def \bx{{\bf x}}
\def \f{{\bf f}}
\def \bX{{\bf X}}
\def \bH{{\bf H}}
\def \bS{{\bf S}}
\def \bF{{\bf F}}
\def \bN{{\bf N}}
\def \bK{{\bf K}}
\def \bE{{\bf E}}
\def \bB{{\bf B}}
\def \bQ{{\bf Q}}
\def \bd{{\bf d}}
\def \bY{{\bf Y}}
\def \bU{{\bf U}}
\def \bL{{\bf L}}
\def \bQ{{\bf Q}}
\def \bP{{\bf P}}
\def \bR{{\bf R}}
\def \bC{{\bf C}}
\def \bD{{\bf D}}
\def \bM{{\bf M}}
\def \bZ{{\mathbb Z}}
\def \xtoleqr{x^{(\leq r)}}
\def \hU{\hat{U}}
\def \k{\kappa}
\def \ee{\overline{p^{\k}}}

\newtheorem{THM}{{\!}}[section]
\newtheorem{THMX}{{\!}}
\renewcommand{\theTHMX}{}
\newtheorem{theorem}{Theorem}[section]
\newtheorem{corollary}[theorem]{Corollary}
\newtheorem{lemma}[theorem]{Lemma}
\newtheorem{proposition}[theorem]{Proposition}
\theoremstyle{definition}
\newtheorem{definition}[theorem]{Definition}
\theoremstyle{remark}
\newtheorem{remark}[theorem]{Remark}
\newtheorem{example}[theorem]{\bf Example}
\numberwithin{equation}{section}
\address{University of New Mexico \\ Albuquerque, NM 87131}
\email{buium@math.unm.edu, arnab@math.unm.edu} 
\maketitle

\begin{abstract}
We prove that some of the basic differential functions appearing in the (unramified) theory of arithmetic differential equations \cite{book}, especially some of the basic differential modular forms in that theory,  arise from a ``ramified situation". This property can be viewed as a special kind of overconvergence property.
One can also go in the opposite direction by using differential functions that arise in  a ramified situation  to construct ``new''   (unramified) differential functions.
\end{abstract}

\section{Introduction}

This paper is a continuation of the study of arithmetic differential equations begun in \cite{char, difmod}; cf. the Introduction and bibliography of \cite{book}. For the convenience of the reader the present paper is written so as to be logically independent of the above references; we will instead quickly review here the main concepts of that theory and we will only refer to \cite{char,difmod,book} for various results as need. The purpose of the theory in \cite{char,difmod,book} is to develop an arithmetic analog of ordinary differential equations. This theory has  a series of purely arithmetic applications for which we refer to  \cite{pjets,difmod,BP}.
In the rest of the introduction we will define our main concepts and  state (in a rough form) our main results. We shall refer to the main body of the paper for detailed statements and for the proofs of our results.

\subsection{Review of notation and terminology \cite{char,difmod,book}}
Throughout  this paper $p\geq 5$ is a fixed prime and we denote by $R_p=\widehat{\mathbb Z}_p^{ur}$ the completion of the maximum unramified extension of ${\mathbb Z}_p$.
We set $K_p=R_p[1/p]$ (fraction field of $R_p$) and $k=R_p/pR_p$ (residue field of $R_p$); so $k$ is  an algebraic closure of ${\mathbb F}_p$.
Let $\pi$ be a root of an Eisenstein polynomial of degree $e\geq 2$ with coefficients in ${\mathbb Z}_p$. (Recall that ${\mathbb Q}_p(\pi)$ is then a totally ramified extension of ${\mathbb Q}_p$; conversely  any finite totally ramified extension of ${\mathbb Q}_p$
 is of the form ${\mathbb Q}_p(\pi)$ with $\pi$ a root of an Eisenstein polynomial with coefficients in ${\mathbb Z}_p$.)
 In order to simplify some of our exposition we will assume in what follows that ${\mathbb Q}_p(\pi)/{\mathbb Q}_p$ is a Galois.
(A typical example we have in mind for our applications is $\pi=1-\zeta_p$ where $\zeta_p$ always denotes in this paper a $p$-th root of unity; in this case  $e=p-1$.)
Consider the ring
$R_{\pi}:=R_p[\pi]=R_p\otimes_{{\mathbb Z}_p}{\mathbb Z}_p[\pi]$. Then $R_{\pi}$ is a complete discrete valuation ring with maximal ideal generated by $\pi$ and with fraction field $K_{\pi}$ of degree $e$ over $K_p$. Let $v_p$ be the $p$-adic valuation on an algebraic closure of $K_{\pi}$ such that
 $v_p(p)=1$  (so $v_p(\pi)=1/e$) and let $|x|:=p^{-v_p(x)}$ be the corresponding absolute value. The ring $R_{\pi}$  possesses a unique ring automorphism $\phi$ such that $\phi(\pi)=\pi$ and $\phi$ lifts the $p$-power Frobenius of $k=R_{\pi}/\pi R_{\pi}$. Clearly $\phi$ sends $R_p$ into itself and is a lift of the $p$-power Frobenius of $k=R_p/pR_p$.
Also throughout the paper $\h$ denotes $p$-adic  completion. For $R_{\pi}$-algebras the $p$-adic completion $\h$ is, of course the same as the $\pi$-adic completion.

Our substitutes for ``differentiation" with respect to $p$ and $\pi$ respectively are the {\it Fermat quotient maps} \cite{char} $\d_p:R_p\ra R_p$ and $\d_{\pi}:R_{\pi}\ra R_{\pi}$ defined by
$$\begin{array}{rcll}
\d_p x & := & \frac{\phi(x)-x^p}{p}, & x\in R_p,\\
\  & \  & \  & \  \\
\d_{\pi} x & := & \frac{\phi(x)-x^p}{\pi}, &  x\in R_{\pi},\end{array}$$
respectively.
In particular, for $x \in R_p$, we have
$$\begin{array}{rcl}
\d_{\pi}x & = & \frac{p}{\pi}\d_p x,\\
\  & \  & \  \\
\d_{\pi}^2 x & = &
\frac{p^2}{\pi^2}\d_p^2 x + \left( \frac{p}{\pi^2}-\frac{p^p}{\pi^{p+1}}\right)
(\d_p x)^p, ...\end{array}$$

Let $V$ be an affine smooth scheme over $R_p$ and fix a closed embedding $V \subset {\mathbb A}^d$ over $R_p$. (The concepts below will be independent of the embedding.) A function $f_p:V(R_p)\ra R_p$ is called a $\d_p$-{\it function} (or order $r\geq 0$) if there exists a restricted power series $F_p$ with $R_p$-coefficients, in $(r+1)d$ variables such that
\begin{equation}
\label{lunchch}f_p(x)=F_p(x,\d_p x,...,\d_p^r x), \ \ x \in V(R_p)\subset R_p^d.\end{equation}
Here and later a power series is called {\it restricted} if its coefficients tend to $0$. (If $V$ is not necessarily affine $f_p$ is called a $\d_p$-function if its restriction to the $R_p$-points of any affine subset of $V$ is a $\d_p$-function.)
A function $f_{\pi}:V(R_{\pi})\ra R_{\pi}$ is called a $\d_{\pi}$-{\it function} (or order $r\geq 0$) if there exists a restricted power series $F_{\pi}$ with $R_{\pi}$-coefficients, in $(r+1)d$ variables such that
\begin{equation}
f_{\pi}(x)=F_{\pi}(x,\d_{\pi} x,...,\d_{\pi}^r x), \ \ x \in V(R_{\pi})\subset R_{\pi}^d.\end{equation}
(If $V$ is not necessarily affine $f_{\pi}$ is called a $\d_{\pi}$-function if its restriction to the $R_{\pi}$-points of any affine subset of $V$ is a $\d_{\pi}$-function.)

\subsection{$\d_{\pi}$-overconvergence}
The main concept we would like to explore (and exploit) in this paper is given in the following definition.
Let $f_p:V(R_p)\ra R_p$ be a $\d_p$-function. We will say that $f_p$ is $\d_{\pi}$-{\it overconvergent} if one of the following equivalent conditions is satisfied:

\medskip

1) There exists an integer $\nu\geq 0$ and a $\d_{\pi}$-function $f_{\pi}$ making the diagram below commutative:
\begin{equation}
\label{pantof1}\begin{array}{rcl}
V(R_p) & \stackrel{p^{\nu}f_p}{\longrightarrow} & R_p\\
\iota \downarrow & \  & \uparrow Tr\\
V(R_{\pi}) & \stackrel{f_{\pi}}{\longrightarrow} & R_{\pi}\end{array}\end{equation}
(Here $\iota$ stands for the inclusion and $Tr$ stands for the $R_{\pi}/R_p$-trace.)

\medskip

 2) There exists an integer $\nu\geq 0$ and a (necessarily unique) $\d_{\pi}$-function $f_{\pi}$ making the diagram below commutative:
\begin{equation}
\label{pantof2}\begin{array}{rcl}
V(R_p) & \stackrel{p^{\nu}f_p}{\longrightarrow} & R_p\\
\iota \downarrow & \  & \downarrow \iota\\
V(R_{\pi}) & \stackrel{f_{\pi}}{\longrightarrow} & R_{\pi}\end{array}\end{equation}

\medskip

The equivalence between conditions 1 and 2 above is trivial to check; cf. also Proposition
\ref{cconverse}.

The concept of $\d_{\pi}$-overconvergence is related to the classical concept of {\it overconvergence} in the theory of Dwork, Monsky and Washnitzer. Indeed
 let us say that a $\d_p$-function $f_p:V(R_p)\ra R_p$  as in (\ref{lunchch}) is {\it $\d_p$-overconvergent with radius $\geq \rho$}
if for any affine cover of $V$ and any affine embeddings of the open sets of the cover the series $F_p$ in (\ref{lunchch})
can be chosen to be  overconvergent (in the classical sense of Dwork, Monsky and Washnitzer) in the variables $\d_p x,...,\d_p^r x$ ``with radius $\geq \rho$". See the body of the paper for details of this definition. We will then show that any
$\d_p$-function of order $r \leq e-1$ which is $\d_{\pi}$-overconvergent must be $\d_p$-overconvergent with radius greater than or equal to  a universal constant that depends only on $p$ and $e$.

\subsection{Main results}
The interaction between
 $\d_p$-functions and $\d_{\pi}$-functions turns out to be  a two way avenue as follows:

\medskip

 1) {\it From $\d_{\pi}$-functions to $\d_p$-functions}. Given a $\d_{\pi}$-function $f_{\pi}:V(R_{\pi})\ra R_{\pi}$ the function $f_p$ defined by the diagram (\ref{pantof1}) with $\nu=0$ turns out to be a $\d_p$-function.
 In this paper we will
  construct   ``interesting'' $\d_{\pi}$-functions
using bad reduction phenomena and then we will apply
trace constructions (a {\it geometric trace} construction and also
the $R_{\pi}/R_p$-trace construction in diagram (\ref{pantof1}) which can be referred to as an {\it arithmetic trace})  to get ``new"  $\d_p$-functions. Cf. Theorem \ref{walk1}.

\medskip

2) {\it From $\d_p$-functions to $\d_{\pi}$-functions}.
In this paper
we discover that some of the basic ``old" $\d_p$-functions that played a role
in \cite{char,difmod,book}
are  $\d_{\pi}$-overconvergent. Cf. Theorem \ref{walk2}.

\medskip

We will apply the above considerations  mainly to the theory of differential modular forms \cite{difmod, book}. To explain this recall the modular curve $X_1(N)_{R_p}$ over $R_p$ with $(N,p)=1$, $N>4$. This curve is smooth and carries a line bundle $L$ such that the spaces of sections $H^0(X_1(N)_{R_p},L^{\kappa})$ identify with the spaces of modular forms on $\Gamma_1(N)$ defined over $R_p$ of weight $\kappa$; cf. \cite{Gross}, p. 450, where $L$ was denoted by $\omega$.
The curve $X_1(N)_{R_p}$ contains two remarkable (disjoint) closed subsets: the {\it cusp locus}  $(cusps)$ and the {\it supersingular locus} $(ss)$.
On $Y_1(N)=X_1(N)\backslash (cusps)$ the line bundle $L$ identifies with $u_*\Omega^1_{E/Y_1(N)}$  where $u:E\ra Y_1(N)$ is the corresponding universal elliptic curve.
Next consider
  an affine open set $X \subset X_1(N)_{R_p}$  and consider the restriction of $L$ to $X$ which we continue to denote by $L$. We can  consider
the affine $X$-scheme $V:=Spec\left(\bigoplus_{n \in {\mathbb Z}}L^{\otimes n}\right)\ra X$. Then
a $\d_p$-{\it modular function} (on $X$, of level $N$ and order $r$) is simply a $\d_p$-function $V(R_p)\ra R_p$ (of order $r$).
 Similarly a $\d_{\pi}$-{\it modular function} (on $X$, of level $N$ and order $r$) is a $\d_{\pi}$-function $V(R_{\pi})\ra R_{\pi}$ (of order $r$). There is a natural concept of weight for a $\d_p$-modular function or a $\d_{\pi}$-modular function; weights are elements in the ring ${\mathbb Z}[\phi]$ of polynomials in $\phi$ with ${\mathbb Z}$-coefficients; cf. the body of the text for the definition of weight. $\d_p$-modular functions (respectively $\d_{\pi}$-modular functions) possessing weights are called $\d_p$-{\it modular forms} (respectively $\d_{\pi}$-{\it modular forms}).
Now, as we shall review in the body of the paper, $\d_p$-modular functions $f$ (and hence forms) possess $\d_p$-{\it Fourier expansions} denoted by $E(f)$ which are
restricted power series in variables $\d_p q,...,\d_p^r q$, with coefficients in the  ring $R_p((q))\h$.

Our first main result
 is  a construction of some interesting ``new''  $\d_p$-modular forms as $R_{\pi}/R_p$-traces
 of some $\d_{\pi}$-modular forms. In their turn, these $\d_{\pi}$-modular forms will be  constructed using the bad reduction of modular curves. Here is the result (in which $X$ is assumed to be disjoint from  the supersingular locus):

 \begin{theorem}
 \label{walk1}
Let $f=\sum a_n q^n$ be a classical normalized newform  of weight $2$ and level $\Gamma_0(Np)$ over ${\mathbb Z}$. Assume $a_p=1$ and let
$\pi=1-\zeta_p$.
  Then there exists a   $\d_p$-modular form $f^{\sharp}_p$ of level $N$, order $1$, and weight $0$ which is $\d_{\pi}$-overconvergent  and whose $\d_p$-Fourier expansion  satisfies
 the following congruence mod $p$:
$$  E(f^{\sharp}_p)
\equiv   \left(\sum_{(n,p)=1}\frac{a_n}{n}q^n\right)-
\left(\sum_{n \geq 1} a_n q^{np}\right) \frac{\d_p q}{q^p}.
$$
\end{theorem}

Cf. Proposition \ref{patru} in the paper. Note that the condition $a_p=1$ is equivalent to the condition  that the elliptic curve attached to $f$ via the Eichler-Shimura construction have split multiplicative reduction at $p$.
The $\d_p$-modular form $f^{\sharp}_p$ in Theorem \ref{walk1}
 should be viewed as a bad reduction analogue of the $\d_p$-modular forms $f^{\sharp}=f^{\sharp}_p$ of level $N$, order $\leq 2$, and weight $0$
 that were attached in \cite{eigen} to  classical  normalized newforms $f=\sum a_n q^n$ of weight $2$ and level $\Gamma_0(N)$ over ${\mathbb Z}$. For such an $f$ on $\Gamma_0(N)$ that does not have CM (in the sense that the elliptic curve attached to it via the Eichler-Shimura construction does not have CM)  the forms $f^{\sharp}_p$ have order exactly $2$ and  were shown in \cite{BP}  to have   $\d_p$-Fourier expansions satisfying the following congruence mod $p$:
$$E(f_p^{\sharp}) \equiv \left(\sum_{(n,p)=1}\frac{a_n}{n} q^n\right) - a_p \left(\sum_{m\geq 1}a_m q^{mp}\right) \frac{\d_p q}{q^p}+  \left(\sum_{m\geq 1}a_m q^{mp^2}\right) \cdot \left(\frac{\d_p q}{q^p}\right)^p.$$
Similar results are available for $f$ on $\Gamma_0(N)$  having CM; cf. \cite{eigen, BP}. Unlike the forms $f_p^{\sharp}$ for $f$ on $\Gamma_0(Np)$ the forms $f_p^{\sharp}$ for $f$ on $\Gamma_0(N)$ were defined for {\it any $X$} (not necessarily disjoint from the supersingular locus).

Our second main result
is a construction of $\d_{\pi}$-modular forms from certain $\d_p$-modular forms.
Indeed, a key role in the theory in \cite{difmod, Barcau, book} was played by certain $\d_p$-modular forms
denoted by $f^1_p,f_p^2,f_p^3,...$ of weights  $-1-\phi,-1-\phi^2,-1-\phi^3,...$ and by   $\d_p$-modular forms denoted by $f^{\partial}_p$ and $f_{\partial,p}$  of weights $\phi-1$ and $1-\phi$ respectively
(where the former are defined whenever $X$ is disjoint from the cusps while the latter are only defined if $X$ is disjoint from both the cusps and the supersingular locus). Recall that $f^{\partial}_p f_{\partial,p}=1$.
The definition of these forms will be reviewed in the body of the paper.
Our second main result (cf. Theorems \ref{overfr}, \ref{overfpartial}, and \ref{lastover} in the body of the paper) is the following:

\begin{theorem}
\label{walk2}
Assume $v_p(\pi)\geq \frac{1}{p-1}$. Then
the $\d_p$-functions
$f^{\partial}_p,f_{\partial,p}, f_p^1,f_p^2,f_p^3,...$  are  $\d_{\pi}$-overconvergent.
Also $f^{\sharp}_p$ is $\d_{\pi}$-overconvergent for any  classical normalized newform $f$  of weight $2$ and level $\Gamma_0(N)$ over ${\mathbb Z}$.
\end{theorem}

By the way the forms $f_p^1,f^{\partial}_p,f_{\partial,p}$ ``generate" (in a sense explained in \cite{difmod,Barcau,book}) all the so called {\it isogeny covariant} $\d_p$-modular forms (in the sense of loc.cit.). We refer to loc.cit. for the role of these forms in the theory and for the significance of the theory itself (in relation, for instance, to the construction in $\d_p$-{\it geometry} of the quotient of the modular curve by the action of the Hecke correspondences); reviewing this background here would take as too far afield and is not necessary for the understanding of our second main result above.

\subsection{Summary of the main forms}
We end our discussion by summarizing (cf. the table below) the main
$\d_{\pi}$-overconvergent
$\d_p$-modular forms appearing in this paper.
\bigskip

\begin{center}
\begin{tabular}{||l|l|c|c|l||} \hline \hline
\  & \  & \  & \ & \  \\
form   & attached to & order  &  weight & domain $X$ \\
\  & \  & \  & \  & \  \\
\hline
\  & \  & \ & \   & \  \\
$f^r_p$ &  $r \geq 1$ & $r$ & $-1-\phi^r$ &  $X$ disjoint from $(cusps)$ \\
\  & \  & \  & \ & \  \\
\hline
\  & \  & \  & \ & \   \\
$f^{\sharp}_p$ &   $f$ on $\Gamma_0(N)$ &  $1$ or $2$ & $0$ & $X$ arbitrary \\
\  & \  & \  & \ & \   \\
\hline
\  & \  & \  & \ & \  \\
$f^{\sharp}_p$ &  $f$ on $\Gamma_0(Np)$  & $1$ & $0$ & $X$ disjoint from $(ss)$\\
\  & \  & \  & \  & \  \\
\hline
\  & \  & \  & \  & \  \\
$f^{\partial}_p$ & \   & $1$ & $\phi-1$ & $X$ disjoint from $(cusps)$ and $(ss)$\\
\  & \  & \  & \  & \  \\
\hline
\  & \  & \  & \  & \  \\
$f_{\partial,p}$ & \   & $1$ & $1-\phi$ & $X$ disjoint from $(cusps)$ and $(ss)$\\
\  & \  & \  & \  & \  \\
\hline
\end{tabular}
\end{center}

\bigskip

\subsection{Plan of the paper}  We begin, in section 2, by revisiting our main set theoretic concepts above from a scheme theoretic viewpoint; $\d_p$-functions and $\d_{\pi}$-functions will appear as functions
on certain formal schemes called $p$-jet spaces and $\pi$-jet spaces respectively; cf. \cite{char,pjets}.
We shall review some of the properties of the latter and we shall analyze the concept of $\d_{\pi}$-overconvergence in some detail. Section 3 is mainly devoted to  reviewing some basic aspects of modular parameterization and bad reduction of modular curves, following \cite{DI,DR,Gross}; so this section is exclusively concerned with ``non-differential''  matters. In section 4 we go back to  arithmetic differential equations:
we will use modular parameterizations and bad reduction of modular curves to construct certain $\d_{\pi}$-modular forms and eventually the ``new" $\d_p$-modular forms in Theorem \ref{walk1}. In section 5 we prove $\d_{\pi}$-overconvergence of some of the basic $\d_p$-functions of the theory, in particular we prove  Theorem \ref{walk2}.

\subsection{Acknowledgment}
This material is based upon work supported by the National Science
Foundation under Grant No. 0852591 and by the Max Planck Institut
f\"{u}r Mathematik, Bonn. Any opinions, findings, and conclusions or
recommendations expressed in this material are those of the author
and do not necessarily reflect the views of the National Science
Foundation or the Max Planck Institut.

\section{$\d_{\pi}$-overconvergence: definition and general properties}

As expained in the Introduction we begin in this section by presenting $\d_p$-functions and $\d_{\pi}$-functions from a scheme-theoretic viewpoint (which is equivalent to the set-theoretic viewpoint of the Introduction).  The scheme-theoretic viewpoint is less direct than the set-theoretic one  but is the correct viewpoint when it comes to   proofs so will be needed in the sequel. We then introduce the concept (and examine some general properties) of
{\it  $\d_{\pi}$-overconvergence}.

\subsection{$p$-jet spaces and $\pi$-jet spaces \cite{char}}

Let $C_p(X,Y) \in \bZ[X,Y]$ be the polynomial with
integer coefficients \[C_p(X,Y):=\frac{X^p+Y^p-(X+Y)^p}{p}.\] A
$p$-{\it derivation} from a ring $A$ into an $A-$algebra $B$,
$\varphi:A \ra B$, is a map $\d_p:A \ra B$ such that $\d_p(1)=0$ and
\[\begin{array}{rcl}
\d_p(x+y) & = &  \d_p x + \d_p y
+C_p(x,y)\\
\d_p(xy) & = & x^p \cdot \d_p y +y^p \cdot \d_p x
+p \cdot \d_p x \cdot \d_p y,
\end{array}\] for all $x,y \in A$. Given a
$p-$derivation we always denote by $\phi:A \ra B$ the map
$\phi(x)=\varphi(x)^p+p \d_p x$; then $\phi$ is a ring homomorphism. A
$\d_p$-{\it prolongation sequence} is a sequence $S^*=(S^n)_{n \geq 0}$ of  rings $S^n$, $n
\geq 0$, together with ring homomorphisms (still denoted by) $\varphi:S^n \ra
S^{n+1}$ and $p-$derivations $\d_p:S^n \ra S^{n+1}$ such that
$\d_p \circ \varphi=\varphi \circ \d_p$ on $S^n$ for all $n$. We view $S^{n+1}$ as an $S^n-$algebra via $\varphi$. A morphism
of $\d_p$-prolongation sequences, $u^*:S^* \ra \tilde{S}^*$ is a sequence
$u^n:S^n \ra \tilde{S}^n$ of ring homomorphisms such that $\delta_p
\circ u^n=u^{n+1} \circ \d_p$ and $\varphi \circ u^n=u^{n+1} \circ
\varphi$. Let $W$ be the ring of polynomials $\bZ[\phi]$ in the
indeterminate $\phi$.  For $w=\sum a_i \phi^i$  (respectively for $w$ with $a_i \geq 0$), $S^*$ a $\d_p$-prolongation
sequence, and $x \in (S^0)^{\times}$ (respectively $x \in S^0$) we
can consider the element $x^w:=\prod_{i=0}^r \varphi^{r-i}
\phi^i(x)^{a_i} \in (S^r)^{\times}$ (respectively $x^w \in S^r$).

Recall the ring $R_p:=\hat{\bZ}_p^{ur}$, completion of the maximum
unramified extension of the ring of $p$-adic integers $\bZ_p$.
 Then $R_p$ has a unique
$p-$derivation $\d_p:R_p \ra R_p$ given by $$\d_p x=(\phi(x)-x^p)/p,$$ where
$\phi:R_p \ra R_p$ is the unique lift of the $p-$power Frobenius map
on $k=R_p/pR_p$. One can consider the  $\d_p$-prolongation sequence $R_p^*$ where
$R_p^n=R_p$ for all $n$. By a $\d_p$-{\it prolongation sequence over $R_p$} we
understand a prolongation sequence $S^*$ equipped with a morphism
$R_p^* \ra S^*$. From now on all our $\d_p$-prolongation sequences are
assumed to be over $R_p$.

Let now $\pi$ be a root of an Eisentein polynomial with ${\mathbb Z}_p$-coefficients and let
 $C_{\pi}(X,Y) \in {\mathbb Z}_p[\pi][X,Y]$ be the polynomial
 \[C_{\pi}(X,Y):=\frac{X^p+Y^p-(X+Y)^p}{\pi}=\frac{p}{\pi}C_p(X,Y).\]
A $\pi$-{\it derivation} from an ${\mathbb Z}_p[\pi]$-algebra
 $A$ into an $A-$algebra $B$,
$\varphi:A \ra B$, is a map $\d_{\pi}:A \ra B$ such that $\d_{\pi}(1)=0$ and
\[\begin{array}{rcl}
\d_{\pi}(x+y) & = &  \d_{\pi} x + \d_{\pi} y
+C_{\pi}(x,y)\\
\d_{\pi}(xy) & = & x^p \cdot \d_{\pi} y +y^p \cdot \d_{\pi} x
+\pi \cdot \d_{\pi} x \cdot \d_{\pi} y,
\end{array}\] for all $x,y \in A$. Given a
$\pi-$derivation we always denote by $\phi:A \ra B$ the map
$\phi(x)=\varphi(x)^p+\pi \d_{\pi} x$; then $\phi$ is a ring homomorphism. A
$\d_{\pi}$-{\it prolongation sequence} is a sequence $S^*=(S^n)_{n \geq 0}$ of  ${\mathbb Z}_p[\pi]$-algebras $S^n$, $n
\geq 0$, together with ${\mathbb Z}_p[\pi]$-algebra homomorphisms (still denoted by) $\varphi:S^n \ra
S^{n+1}$ and $\pi-$derivations $\d_{\pi}:S^n \ra S^{n+1}$ such that
$\d_{\pi} \circ \varphi=\varphi \circ \d_{\pi}$ on $S^n$ for all $n$.  A morphism
of $\d_{\pi}$-prolongation sequences, $u^*:S^* \ra \tilde{S}^*$ is a sequence
$u^n:S^n \ra \tilde{S}^n$ of ${\mathbb Z}_p[\pi]$-algebra
 homomorphisms such that $\delta_{\pi}
\circ u^n=u^{n+1} \circ \d_{\pi}$ and $\varphi \circ u^n=u^{n+1} \circ
\varphi$. Let $W$ be, again,  the ring of polynomials $\bZ[\phi]$ in the
indeterminate $\phi$.  For $w=\sum a_i \phi^i$  (respectively for $w$ with $a_i \geq 0$), $S^*$ a $\d_{\pi}$-prolongation
sequence, and $x \in (S^0)^{\times}$ (respectively $x \in S^0$) we
can consider the element $x^w:=\prod_{i=0}^r \varphi^{r-i}
\phi^i(x)^{a_i} \in (S^r)^{\times}$ (respectively $x^w \in S^r$).

As above we may consider $R_{\pi}=R_p[\pi]$ and the
$\pi-$derivation $\d_{\pi}:R_{\pi} \ra R_{\pi}$ given by $$\d_{\pi} x=(\phi(x)-x^p)/\pi.$$ One can consider the  $\d_{\pi}$-prolongation sequence $R_{\pi}^*$ where
$R_{\pi}^n=R_{\pi}$ for all $n$. By a $\d_{\pi}$-{\it prolongation sequence over $R_{\pi}$} we
understand a prolongation sequence $S^*$ equipped with a morphism
$R_{\pi}^* \ra S^*$. From now on
all our $\d_{\pi}$-prolongation sequences are
assumed to be over $R_{\pi}$.

We note that if $S^*=(S^n)_{n \geq 0}$ is a $\d_p$-prolongation sequence
 such that each $S^n$ is flat over $R_p$ then the sequence
$S^*\otimes_{R_p} R_{\pi}=(S^n \otimes_{R_p} R_{\pi})_{n \geq 0}$
has  a natural structure of $\d_{\pi}$-prolongation sequence.
Indeed letting $\phi:S^n \ra S^{n+1}$ denote, as usual,
the ring homomorphisms $\phi(x)=x^p+p\d_p x$ one can extend these $\phi$s to ring homomorphisms $\phi:S^n \otimes_{R_p}R_{\pi} \ra S^{n+1}\otimes_{R_p}R_{\pi}$
by the formula $\phi(x \otimes y)=\phi(x) \otimes \phi(y)$ where $\phi:R_{\pi}\ra R_{\pi}$ is given, as usual, by $\phi(y)=y^p+\pi \d_{\pi} y$. Then one can define $\pi$-derivations $\d_{\pi}:S^n \otimes_{R_p}R_{\pi}\ra S^{n+1}\otimes_{R_p}R_{\pi}$
by $\d_{\pi}(z)=(\phi(z)-z^p)/\pi$ for $z \in S^n\otimes_{R_p}R_{\pi}$. With these $\d_{\pi}$s the sequence
$S^*\otimes_{R_p}R_{\pi}$ is a $\d_{\pi}$-prolongation sequence.

For any affine $R_p$-scheme of finite type $X=Spec\ A$ there exists
a (unique) $\d_p$-prolongation sequence, $A^*=(A^n)_{n \geq 0}$, with $A^0=A$ such that for any $\d_p$-prolongation sequence $B^*$
 and any $R_p$-algebra homomorphism $u:A \ra B^0$ there exists a unique morphism of $\d_p$-prolongation sequences $u^*:A^* \ra B^*$ with $u^0=u$.
We define the $p$-{\it jet spaces} $J^n_p(X)$ of $X$ as the formal schemes
$J^n_p(X):=Spf\ \hat{A^n}$. This construction immediately globalizes to the case $X$ is not necessarily affine (such that the construction commutes, in the obvious sense, with open immersions). For $X$  smooth over $R_p$
the ring of $\d_p$-functions $X(R_p) \ra R_p$ naturally identifies with the
ring of global functions $\cO(J^n_p(X))$: under this identification any function $f \in \cO(J^n_p(X))$
gives rise to a $\d_p$-function $X(R_p)\ra R_p$ sending any point $P \in X(R_p)$,
$P:Spec\ R_p\ra X$ into the $R_p$-point of the affine line ${\mathbb A}^1_{R_p}$ defined by
$$Spf\ R_p \stackrel{P^n}{\ra} J_p^n(X) \stackrel{f}{\ra} \hat{\mathbb A}^1_{R_p};$$
here $P^n$ is the morphism induced from $P$ via the universality property of the $p$-jet space.  If $X$ is a group scheme over $R_p$ then
$$f:J^n_p(X)\ra \hat{\mathbb G}_{a,R_p}=\hat{\mathbb A}^1_{R_p}$$
 is a group homomorphism into the additive group of the line if and only if the corresponding map $X(R_p)\ra R_p$ is a group homomorphism; such an $f$ is called a $\d_p$-{\it character} of $X$.

As a prototypical example if $X={\mathbb A}^N_{R_p}=Spec\ R_p[x]$ is the affine space (where $x$ is an $N$-tuple of variables) then $J^n_p(X)=Spf\ R_p[x,\d_p x,...,\d_p^n x]\h$ (where $\d_px,...,\d_p^n x$ are new $N$-tuples of variables).

We will need, in this paper, a slight generalization of the above constructions as follows; cf. \cite{difmod}. First note that the $p$-jet spaces $J^n_p(X)$
only depend on the $p$-adic completion of $X$ and not on $X$. This immediately implies that one can introduce $p$-jet spaces $J^n_p({\mathcal X})$ attached
to formal $p$-adic schemes ${\mathcal X}$ over $R_p$ which are locally $p$-adic completions of schemes of finite type over $R_p$; the latter association is functorial in ${\mathcal X}$.

Similarly, for any affine $R_{\pi}$-scheme of finite type $Y=Spec\ A$ there exists
a (unique) $\d_{\pi}$-prolongation sequence, $A^*=(A^n)_{n \geq 0}$, with $A^0=A$ such that for any $\d_{\pi}$-prolongation sequence $B^*$
 and any $R_{\pi}$-algebra homomorphism $u:A \ra B^0$ there exists a unique morphism of $\d_{\pi}$-prolongation sequences $u^*:A^* \ra B^*$ with $u^0=u$.
We define the $\pi$-{\it jet spaces} $J^n_{\pi}(Y)$ of $Y$ as the formal schemes
$J^n_{\pi}(Y):=Spf\ \hat{A^n}$. This construction immediately globalizes to the case $Y$ is not necessarily affine (such that the construction commutes, in the obvious sense, with open immersions). Again, for $Y$  smooth over $R_{\pi}$
the ring of $\d_{\pi}$-functions $Y(R_{\pi}) \ra R_p$ naturally identifies with the
ring of global functions $\cO(J^n_{\pi}(Y))$.
 If $Y$ is a group scheme over $R_{\pi}$ then $f:J^n_{\pi}(Y)\ra \hat{\mathbb A}^1_{R_{\pi}}$ is a group homomorphism into the additive group of the line if and only if the corresponding map $Y(R_{\pi})\ra R_{\pi}$ is a group homomorphism; such an $f$ is called a $\d_{\pi}$-{\it character} of $Y$.

As a prototypical example if $Y={\mathbb A}^N_{R_{\pi}}=Spec\ R_{\pi}[x]$ is the affine space  then $J^n_{\pi}(Y)=Spf\ R_{\pi}[x,\d_{\pi} x,...,\d_{\pi}^n x]\h$ (where $\d_{\pi}x,...,\d_{\pi}^n x$ are new $N$-tuples of variables).

As in the case of $p$-jet spaces,  note that the $\pi$-jet spaces $J^n_{\pi}(Y)$
only depend on the $\pi$-adic completion of $Y$ and not on $Y$. This immediately implies that one can introduce $\pi$-jet spaces $J^n_{\pi}({\mathcal Y})$ attached
to formal $\pi$-adic schemes ${\mathcal Y}$ over $R_{\pi}$ which are locally $\pi$-adic completions of schemes of finite type over $R_{\pi}$; the latter association is functorial in ${\mathcal Y}$.

 For any scheme $X/R_p$ we write $X_{R_{\pi}}:=X \otimes_{R_p} R_{\pi}$. Let $X/R_p$ be a smooth affine scheme. The $\d_p$-prolongation sequence
 $(\cO(J^n(X)))_{n \geq 0}$
induces a structure of $\d_{\pi}$-prolongation sequence on
the sequence  $(\cO(J^n(X)) \otimes_{R_p} R_{\pi})_{n \geq 0}$.
By the universality property of the $\d_{\pi}$-prolongation sequence
 $(\cO(J^n_{\pi}(X_{R_{\pi}})))_{n \geq 0}$ we get a  canonical morphism of $\d_{\pi}$-prolongation sequences
\begin{equation}
\label{tra}\cO(J^n_{\pi}(X_{R_{\pi}}))\ra \cO(J^n_p(X))\otimes_{R_p} R_{\pi}.\end{equation}
The following  is trivial to prove by induction:

\begin{lemma} For any $n \geq 1$ there exists
a polynomial $F_n\in R_{\pi}[t_1,...,t_n]$ without constant term, of degree $\leq p^{n-1}$
with the property that for any
  $f \in \cO(X)$ we have
\begin{equation}
\label{macac}
\d_{\pi}^n f \mapsto \frac{p^n}{\pi^n} \d_p^n f +\pi^{\max\{e-n,0\}}F_n(\d_p f,...,\d_p^{n-1}f)
\end{equation}
under the map (\ref{tra}).\end{lemma}

 In particular, for instance,
\begin{equation}
\label{formul}
\begin{array}{rcl}
\d_{\pi} f& \mapsto & \frac{p}{\pi}\d_p f\\
\  & \  & \  \\
\d_{\pi}^2 f & \mapsto & \frac{p^2}{\pi^2}\d_p^2 f + \left( \frac{p}{\pi^2}-\frac{p^p}{\pi^{p+1}}\right)
(\d_p f)^p.\end{array}\end{equation}
Note that for $1 \leq n \leq e-1$ and $f \in \cO(X)$ the image
 of $\d_{\pi}^n f$
in  $\cO(J^n_p(X))\otimes_{R_p} R_{\pi}$ is  always in the ideal generated by $\pi$.
Also note that for $f \in \cO(X)$ the image of $\d_{\pi}^e f$
in  $\cO(J^n_p(X))\otimes_{R_p} R_{\pi}$ is {\it not} always in the ideal generated by $\pi$; indeed the image of $\d^e_{\pi} p$ in $\cO(J^n_p(X))\otimes_{R_p} R_{\pi}$ belongs to $R_{\pi}^{\times}$.
For $X$ not necessarily affine we get a morphism (\ref{tra})  and a canonical morphism of $\pi$-adic formal schemes
\begin{equation}
J^n_p(X)\otimes_{R_p} R_{\pi} \ra J^n_{\pi}(X_{R_{\pi}}).\end{equation}
Note that the map (\ref{tra}) is an isomorphism if $n=0$.
For $n \geq 1$ the map (\ref{tra}) is not  surjective
and its reduction mod $p$ is not injective.
Nevertheless, we have:

\begin{proposition}
\label{inj}
The map (\ref{tra}) is injective.
\end{proposition}

We will usually view the map (\ref{tra}) as an inclusion.

\bigskip

{\it Proof}. Indeed it is enough to prove this for $X$ affine and sufficiently small. So let us assume that $X$ has \'{e}tale coordinates i.e. there is an \'{e}tale map $R[x]\ra \cO(X)$
with $x$ a tuple of variables. Then by the local product property of $\pi$-jet spaces
\cite{char}, Proposition 1.4, (\ref{tra}) becomes the natural map
\begin{equation}
\label{tralala}
\cO(X_{R_{\pi}})[\d_{\pi} x,...,\d_{\pi}^n x]\h \ra \cO(X)[\d_p x,...,\d_p^n x]\h \otimes_{R_p} R_{\pi}.\end{equation}
Now let  $L$ be the fraction field of the $\pi$-adic completion
of $\cO(X_{R_{\pi}})$. (The latter is an integral domain by the smoothness of $X/R_p$.)
Then the left hand side of (\ref{tralala}) embeds into
$L[[\d_{\pi} x,...,\d_{\pi}^n x]]$ while
the right hand side of (\ref{tralala}) embeds into
$L[[\d_p x,...,\d_p^n x]]$
(the latter because $R_{\pi}$ is finite over $R_p$).
Finally we claim that we have a natural isomorphism
\begin{equation}
\label{ddda}
L[[\d_{\pi} x,...,\d_{\pi}^n x]]\simeq L[[\d_p x,...,\d_p^n x]]\end{equation}
that induces (\ref{tra}); this of course will end the proof that (\ref{tra}) is injective. To prove the claim note that there is  natural homomorphism
\begin{equation}
\label{tworings}
L[\d_{\pi} x,...,\d_{\pi}^n x]\ra L[\d_p x,...,\d_p^n x]\end{equation}
which is trivially seen by induction to be surjective by the formulae (\ref{macac}). Since the rings in (\ref{tworings}) have both dimension $n$ it follows that (\ref{tworings}) is an isomorphism. Since (\ref{tworings}) maps the ideal
$(\d_{\pi} x,...,\d_{\pi}^n x)$ into (and hence onto) the ideal $(\d_p x,...,\d_p^n x)$ we get an isomorphism like in (\ref{ddda}) and we are done.
\qed

\medskip

Let now $Tr:R_{\pi}\ra R_p$ be the ($R_p$-linear) trace map. We may consider the $R_p$-linear map
\begin{equation}
\label{Tr}
1 \otimes Tr: \cO(J^n_p(X))\otimes_{R_p} R_{\pi}\ra \cO(J^n_p(X))\otimes_{R_p} R_p=\cO(J^n_p(X)).
\end{equation}
Composing (\ref{tra}) with (\ref{Tr}) we get an $R_p$-linear {\it arithmetic trace} map:
\begin{equation}
\label{tr}
\tau_{\pi}:\cO(J^n_{\pi}(X_{R_{\pi}}))\ra \cO(J^n_p(X)).
\end{equation}
(Later we will encounter another type of trace maps which will be referred to as {\it geometric trace} maps.)

\begin{proposition}
\label{cconverse}
Let $X$ be a smooth scheme over $R_p$ and $f \in \cO(J^n_p(X))$. The following conditions are equivalent:

1) $f$ times a power of $p$ belongs to the image of the trace map (\ref{tr}).

2) $f$ times a power of $p$ belongs to the image of the inclusion map (\ref{tra}).
\end{proposition}

{\it Proof}.
The fact that condition 2 implies condition 1 is trivial.

In order to check that condition 1 implies condition 2 let $\Sigma$ be the Galois group of ${\mathbb Q}_p(\pi)/{\mathbb Q}_p$ (and hence also of $K_{\pi}/K_p$) and let us consider the action of $\Sigma$ on $\cO(J^r_p(X))\otimes_{R_p}R_{\pi}$
via the action on the second factor. We will prove that
   $\Sigma$
acts on  the image of  $\cO(J^r_{\pi}(X_{R_{\pi}})) \ra \cO(J^r_p(X))\otimes_{R_p}R_{\pi}$; this will of course end the proof of the Proposition.
Let $\sigma \in \Sigma$ and
$\frac{\sigma \pi}{\pi}=:u\in {\mathbb Z}_p[\pi]^{\times}$.

\medskip

{\it Claim 1.}  $\phi$ and $\sigma$ commute on $R_{\pi}$.
Indeed $\phi \circ \sigma$ and $\sigma \circ \phi$ have the same effect on $R_p$ and on $\pi$.

\medskip

{\it Claim 2}.  $\phi \circ \sigma=\sigma\circ \phi$ as
 maps from $\cO(J^i_p(X))\otimes_{R_p}R_{\pi}$ to
$\cO(J^{i+1}_p(X))\otimes_{R_p}R_{\pi}$.
Indeed it is enough to check this for $X=Spec\ R_p[x]$ the affine space, $x$ a tuple of variables.
So it is enough to check that $\phi$ and $\sigma$ commute as maps
from $R_{\pi}[x,\d_p x,...,\d_p^i x]\h$ to $R_{\pi}[x,\d_p x,...,\d_p^{i+1} x]\h$. This is clear because $\phi$ and $\sigma$ commute on $R_{\pi}$ and on each tuple $\d_p^j x$.

\medskip

{\it Claim 3}. $\sigma \circ \d_{\pi}=\frac{1}{u}\cdot  \d_{\pi} \circ \sigma$
 as maps from $\cO(J^i_p(X))\otimes_{R_p}R_{\pi}$ to
$\cO(J^{i+1}_p(X))\otimes_{R_p}R_{\pi}$.
This follows trivially from Claim 2.

\medskip

Now to conclude
it is enough to show that for  any $1\leq i \leq r$, and any $f \in \cO(X)$ we have that $\sigma (\d_{\pi}^i f)$
is obtained by evaluating a polynomial $P_i$
with $R_{\pi}$-coefficients at $\d_{\pi}f,...,\d_{\pi}^i f$.
We proceed by induction on $i$. The case $i=1$ is clear. Assume our assertion is true for $i$.
Then
$$\begin{array}{rcl}
\sigma \d_{\pi}^{i+1} f & = & \sigma \d_{\pi} (\d_{\pi}^i f)\\
\  & \  & \  \\
\ & = & \frac{1}{u} \d_{\pi} (\sigma(\d_{\pi}^i f))\\
\  & \  & \  \\
\  & = & \frac{1}{u} \d_{\pi}(P_i(\d_{\pi}f,...,\d_{\pi}^i f))\end{array}$$
and we are done.
\qed

\begin{definition}
\label{basicdef}
A function $f \in \cO(J^n_p(X))$ is called $\d_{\pi}$-{\it overconvergent} if
it satisfies
one of the equivalent conditions in Proposition \ref{cconverse}.
\end{definition}

\begin{remark}
\label{derr}
The set of $\d_{\pi}$-overconvergent elements of $\cO(J^n_p(X))$ is a subring containing all the elements of the form $\d^i_p f$ with $i \leq n$ and $f \in \cO(X)$. In particular if $X$ is affine then the subring of $\d_{\pi}$-overconvergent elements of $\cO(J^n_p(X))$ is $p$-adically dense in $\cO(J^n_p(X))$.
\end{remark}

\begin{remark}
Under the identification of $\d_p$-functions (respectively $\d_{\pi}$-functions) with elements of the ring $\cO(J^n_p(X))$ (respectively $\cO(J^n_{\pi}(X_{R_{\pi}}))$) the definition of $\d_{\pi}$-overconvergence above corresponds to the definition of $\d_{\pi}$-overconvergence given in the Introduction.
\end{remark}

\begin{remark}
\label{usefull}
Let us note that $\d_{\pi}$-overconvergence is preserved by precomposition with regular maps. Indeed, if $u:Y \ra X$ is a morphism of smooth $R_p$-schemes and if $\lambda \cdot f$ is in the image of (\ref{tra})
for some $\lambda \in R_{\pi}$ and some $f \in \cO(J^r_p(X))$ then
if $f$
is identified with the corresponding map $f:X(R_p)\ra R_p$ it follows that
 $\lambda \cdot f \circ u$ is in the image of
$\cO(J^n_{\pi}(Y_{R_{\pi}}))\ra \cO(J^n_p(Y))\otimes_{R_p} R_{\pi}.$
(Here $f\circ u$ is identified with  $u^* f$ where $u^*$ is the naturally induced map $\cO(J^r_p(X))\ra \cO(J^r_p(Y))$.
\end{remark}

The next Proposition shows that the trace map $\tau_{\pi}$ in (\ref{tra}), although not injective, is ``as close as possible" to being so.

\begin{proposition}
 The map
\begin{equation}
\label{maha}
\cO(J^n_{\pi}(X_{R_{\pi}}))\ra \bigoplus_{i=0}^{e-1} \cO(J^n_p(X)),\ \ \ f \mapsto (\tau_{\pi}(f),\tau_{\pi}(\pi f),...,\tau_{\pi}(\pi^{e-1} f))\end{equation}
is injective.
\end{proposition}

{\it Proof}. Indeed if the image of $f$ in $\cO(J^n_p(X))\otimes_{R_p} R_{\pi}$ is
$\sum_{i=0}^{e-1} f_i \otimes \pi^i$ and the image of $f$ via the map (\ref{maha}) is $0$ then
we get $\sum_{i=0}^{e-1} Tr(\pi^{i+j})f_i=0$ for all $j=0,...,e-1$. Now
$det(Tr(\pi^{i+j}) \neq 0$
which implies $f_0=...=f_{e-1}=0$ hence, by Proposition \ref{inj}, $f=0$.
\qed

\begin{example}
Consider the multiplicative group
${\mathbb G}_{m,R_{p}}=Spec\ R_p[x,x^{-1}]$ and the
standard $\d_p$-character  $$\psi_p \in \cO(J^1_p({\mathbb G}_{m,R_p}))=R_p[x,x^{-1},\d_p x]\h$$ in \cite{char}
defined by
$$\psi_p:=``\frac{1}{p}\log\left(\frac{\phi(x)}{x^p}\right)":=
\sum_{n \geq 1} (-1)^{n-1}\frac{p^{n-1}}{n}\left(\frac{\d_p x}{x^p}\right)^n.$$
 Assume  $v_p(\pi)\geq \frac{1}{p-1}$, e.g. $\pi=1-\zeta_p$. Then clearly $p\psi_p=\pi \psi_{\pi}$ where $\psi_{\pi}$ is the
$\d_{\pi}$-character $\psi_{\pi}\in \cO(J^1_{\pi}({\mathbb G}_{m,R_{\pi}}))$ defined by
\begin{equation}
\label{sol}
\psi_{\pi}:=
\sum_{n \geq 1} (-1)^{n-1}\frac{\pi^{n-1}}{n}\left(\frac{\d_{\pi} x}{x^p}\right)^n.
\end{equation}
(which is well defined because if $v_p(\pi)\geq \frac{1}{p-1}$ then $v_p(\pi^{n-1}/n)$ is $\geq 0$ and $\ra \infty$ as $n \ra \infty$).
 So $\psi_p$ is $\d_{\pi}$-overconvergent. Moreover
$$\tau_{\pi}(\psi_{\pi})= Tr\left(\frac{1}{\pi}\right)\cdot p\psi_p.$$
By the way, if $\pi=1-\zeta_p$ then $Tr(\frac{1}{\pi})=\frac{p-1}{2}$.
\end{example}

The above global concepts and remarks have a local counterpart as follows. Let $q$ be a variable
and $\d_{\pi}^i q$, $\d_p^i q$ corresponding variables. Then exactly as above we have that
the natural map
\begin{equation}
\label{tralalu}
 R_{\pi}((q))[\d_{\pi} q,...,\d_{\pi}^n q]\h \ra R_{\pi}((q))[\d_p q,...,\d_p^n q]\h\end{equation}
is injective. We shall view this map as an inclusion.
On the other hand there is a natural {\it trace map}
 \begin{equation}
 \label{vopp}
 \tau_{\pi}:R_{\pi}((q))[\d_{\pi} q,...,\d_{\pi}^n q]\h \ra
 R_{\pi}((q))[\d_p q,...,\d_p^n q]\h \stackrel{Tr}{\ra}
 R_p((q))[\d_p q,...,\d_p^n q]\h,
 \end{equation}
 where the first map is the inclusion (\ref{tralalu})
  and the second map $Tr$ is induced by the trace map $Tr:R_{\pi}\ra R_p$ on the coefficients of the series. As in the global case we have that:

\begin{proposition}
\label{miedor}
For a series $f$ in $R_p((q))[\d_p q,...,\d^n_p q]\h$ the following conditions are equivalent:

1) $f$ times a power of $p$ belongs to the image of the trace map (\ref{vopp}).

2) $f$ times a power of $p$ belongs to the image of the inclusion map (\ref{tralalu}).
\end{proposition}

So as in the global case we can make the following:

\begin{definition}
A series in
 $R_p((q))[\d_p q,...,\d^n_p q]\h$
is $\d_{\pi}$-{\it overconvergent} if it satisfies one of the equivalent conditions in Proposition \ref{miedor}.
\end{definition}

\begin{example}
Assume $v_p(\pi)\geq \frac{1}{p-1}$, e.g. $\pi=1-\zeta_p$.
Then the series
\begin{equation}
\label{eqqq}
\Psi_p:=``\frac{1}{p}\log \left( \frac{\phi(q)}{q^p}\right)'':=
\sum_{n \geq 1} (-1)^{n-1}\frac{p^{n-1}}{n} \left(\frac{\d_p q}{q^p}\right)^n
\in R_p((q))[\d_p q]\h\end{equation}
is $\d_{\pi}$-overconvergent. Indeed we can write
$p\Psi_p=\pi\Psi_{\pi}$ where the series
$$\Psi_{\pi}:=
\sum_{n \geq 1} (-1)^{n-1}\frac{\pi^{n-1}}{n} \left(\frac{\d_{\pi} q}{q^p}\right)^n;$$
is in $R_{\pi}((q))]\d_{\pi}q]\h$ because $v_p(\pi^{n-1}/n)$ is $\geq 0$ and $\ra \infty$.
\end{example}

Next we would like to compare the concept of $\d_{\pi}$-overconvergence introduced above with the classical concept of overconvergence as it was introduced in the work of Dwork, Monsky, and Washnitzer.  Let us recall the classical concept of {\it overconvergence} of power series or, more generally the concept of {\it overconvergence of series with respect to a subset of variables}.

\begin{definition}
Let $C$ be a positive real number and $\rho=p^C$.
Let $x$ and $y$ be tuples of variables and $F \in R_p[x,y]\h\subset R_p[[x,y]]$ a restricted power series, $F=\sum a_{\alpha,\beta} x^{\alpha}y^{\beta}$ (where $\alpha, \beta$ are multiindices). Then $F$ is called {\it overconvergent in the variables $y$ with radius $\geq \rho$}
if there exists a positive real number $C'$ such that for all $\alpha, \beta$ one has
$v_p(a_{\alpha,\beta}) \geq C|\beta|-C'$ (equivalently $|a_{\alpha,\beta}|\rho^{|\beta|}$ is bounded from above independently of $\alpha$ and $\beta$).
\end{definition}

Here $|\beta|$ is, of course, the sum of the components of $\beta$.
We make then the following:

\begin{definition}
\label{odefmica}
 Let $V$ be any smooth scheme over $R_p$ and let $f \in \cO(J^r_p(V))$. We say that $f$ is  {\it $\d_p$-overconvergent with radius $\geq \rho$} if for any affine open set $W\subset V$  and for any closed embedding $W \subset
  {\mathbb A}^d=Spec\ R_p[x]$  (where $x$ is a $d$-tuple of variables)
  the image of $f$ in $\cO(J^r_p(W))$ is  the image via $\cO(J^r_p({\mathbb A}^d)) \ra \cO(J^r_p(V))$ of a restricted power series in $\cO(J^r_p({\mathbb A}^d))=R_p[x,\d_p x,...,\d_p^r x]\h$  which is overconvergent in the variables $\d_p x,...,\d_p^r x$ with radius $\geq \rho$.
\end{definition}

Then we have the following:

\begin{proposition}
\label{caff}
Let $V$ be any smooth $R_p$-scheme.
Assume $1 \leq r \leq e-1$ and $f \in \cO(J^r_p(V))$ is $\d_{\pi}$-overconvergent.
Then $f$ is $\d_p$-overconvergent with radius $\geq p^{(p^{r-1}e)^{-1}}.$
\end{proposition}

{\it Proof}.
It is enough to prove this for $V={\mathbb A}^d$.
By hypothesis  $p^{\nu}F$ is in the image
of
$$R_{\pi}[x,\d_{\pi}x,...,\d_{\pi}^r x]\h \ra R_p[x,\d_p x,...,\d_p^r x]\h
\otimes_{R_p}R_{\pi}$$
for some $\nu$. We may assume $\nu=0$.
 Write
$$F(x,\d_p x,...,\d_p^r x)=\sum_{\alpha_0,...,\alpha_r}a_{\alpha_0...\alpha_r}x^{\alpha_0}
(\d_{\pi}x)^{\alpha_1}...(\d_{\pi}^r x)^{\alpha_r},$$
with $a_{\alpha_0...\alpha_r}\in R_{\pi}$. By
(\ref{macac})
one can find polynomials $G_i\in R_{\pi}[t_1,...,t_i]$ of degree $\leq p^{i-1}$ such that
$$\d_{\pi}^i x=\pi\cdot G_i(\d_p x,...,\d_p^i x),\ \ \ 1\leq i \leq e-1.$$
We get
that
$$F(x,\d_p x,...,\d_p^r x)=\sum a_{\alpha_0...\alpha_r} \pi^{|\alpha_1+...+\alpha_r|} x^{\alpha_0}(G_1(\d_p x))^{\alpha_1}...(G_n(\d_p x,...,\d_p^r x))^{\alpha_r}.$$
Then clearly the coefficient of the monomial
$$x^{\alpha}(\d_p x)^{\beta_1}...(\d_p^r x)^{\beta_r}$$
 in $F$ is going to have $p$-adic valuation at least
$$\frac{1}{e}\left(\frac{|\beta_1+...+\beta_r|}{p^{r-1}}-1\right)$$
and we are done.
\qed

\begin{remark}
Proposition \ref{caff} fails if we do not asssume the order $r$ is strictly less than the ramification index $e$. Here is a typical example. Let $V={\mathbb A}^1=Spec\ R_p[x]$, $\pi=\sqrt{p}$; so $e=2$ and $\d_{\pi}^2 x=p(\d_p^2 x)+u(\d_p x)^p$, $u=1-p^{(p-1)/2}$.
Let $a_n \in R_p$, $v_p(a_n)\ra \infty$, $v_p(a_n)\leq n^{\epsilon}$, $0 < \epsilon <1$,  and let
$$F=F(x,\d_p x, \d_p^2 x):=\sum a_n (p(\d_p^2 x)+u(\d_p x)^p)^n\in R_p[x,\d_p x,\d_p^2 x]\h.$$
Then  $F$ is $\d_{\pi}$-overconvergent because
$$F=\sum a_n (\d_{\pi}^2 x)^n \in R_{\pi}[x,\d_{\pi}x, \d_{\pi}^2 x]\h.$$
On the other hand $F(x,\d_p x,\d_p^2 x)$ is not $\d_p$-overconvergent of radius $\rho$ (regardless of the value of $\rho$). Indeed if this were the case then
$$F(0,y,0)=\sum a_n u^n y^{np}$$
would be  overconvergent in the variable $y$ with radius $\geq \rho$ which is clearly not the case.
\end{remark}

\begin{remark}
The concept of $\d_p$-overconvergence introduced above comes with  a built in  number $\rho$ that is ``coordinate independent" (independent of the affine embedding). It is worth remarking that no such coordinate independent  $\rho$ can be similarly  attached to the classical overconvergence concept of Monski and Washnitzer. This is best exemplified by the following example.
Let $y=\{y_1,y_2\}$ be a pair of variables. For any $R_p$-automorphism $\sigma$ of $R_p[y_1,y_2]$ let us denote by $\hat{\sigma}$ the induced automorphism of $R_p[y_1,y_2]\h$. Then one can easily find examples
of elements
 $f \in R_p[y_1,y_2]\h$ having the following properties:

 1) There is a real $\rho_0>1$ such that the series $f$ is overconvergent in the variables $y$ with radius $\geq \rho_0$;

 2) There is no real $\rho>1$ such that for any automorphism $\sigma$ of $R_p[y_1,y_2]$ the series $\hat{\sigma}(f)$ is overconvergent in the variables $y$ with radius $\geq \rho$.

 To come up with an explicit example let $f=\sum_{n \geq 1} p^n y_1^n$. Then 1) above is satisfied with $\rho_0=p$. To show that  2) is satisfied assume there is a constant $\rho>1$ having the property that for any automorphism $\sigma$ of $R_p[y_1,y_2]$ the series $\hat{\sigma}(f)$ is overconvergent in the variables $y$ with radius $\geq \rho$ and let us seek a contradiction. Take any integer $m \geq 1$ such that $\rho^m >p$ and  let $\sigma$ be the automorphism defined by $\sigma(y_1)=y_1+y_2^m$, $\sigma(y_2)=y_2$. Then $\hat{\sigma}(f)=\sum_{n \geq 1} p^n(y_1+y_2^m)^n$.
 The coefficient of $y_2^{mn}$ in the latter is $p^n$ so we must have that $p^{-n}\rho^{mn}$
 is bounded from above independently of $n$; this is a contradiction and we are done.
\end{remark}

\begin{remark}
Let us make the following definition. For a smooth affine scheme $V$, an element $f \in \cO(J^r_p(V))$
is {\it overconvergent} if there exists a closed embedding $V \subset {\mathbb A}^d=Spec\ R_p[x]$ and a real $\rho >1$
such that  $f$ is  the image via $\cO(J^r_p({\mathbb A}^d)) \ra \cO(J^r_p(V))$ of a restricted power series in $\cO(J^r_p({\mathbb A}^d))=R_p[x,\d_p x,...,\d_p^r x]\h$  which is overconvergent in the variables $x,\d_p x,...,\d_p^r x$ with radius $\geq \rho$. (So here, as opposed to the definition of $\d_p$-overconvergence,  we include the variables $x$ as well. Also we only ask that $\rho$ works for {\it one} particular embedding $V \subset {\mathbb A}^d$.) If $f$ is as above and $V \subset {\mathbb A}^{\tilde{d}}=Spec\ R_p[\tilde{x}]$ is another closed embedding then one can easily see that there exists a generally different  real number $\tilde{\rho}>1$ such that $f$ is  the image via $\cO(J^r_p({\mathbb A}^{\tilde{d}})) \ra \cO(J^r_p(V))$ of a restricted power series in $\cO(J^r_p({\mathbb A}^{\tilde{d}}))=R_p[\tilde{x},\d_p \tilde{x},...,\d_p^r \tilde{x}]\h$  which is overconvergent in the variables $\tilde{x},\d_p \tilde{x},...,\d_p^r \tilde{x}$ with radius $\geq \tilde{\rho}$.)
We expect that all the remarkable $\d_p$-functions appearing in this paper that will be proved to be $\d_{\pi}$-overconvergent (such as $f^r_p,f^{\partial}_p,f_p^{\sharp}$ in the Introduction) are also overconvergent.
This would not  imply (and would not be implied by) our $\d_{\pi}$-overconvergence results or our $\d_p$-overconvergence results with radius bounded by a universal constant.
\end{remark}

\subsection{$p$-jets and $\pi$-jets  of formal groups}
In what follows we recall from \cite{book}, section 4.4, the construction of $p$-jets of formal groups
and  we also introduce the $\pi$-jet analogue of that construction.

Start with a formal group law ${\mathcal F}\in S[[T_1,T_2]]$ (in one variable $T$) over $S=\cO(X)$, where $X$  is a smooth affine $R_p$-scheme. For $r \geq 1$ we let $S^r_p:=\cO(J^r_p(X))$. Let ${\bf T}$ be the pair of variables $T_1,T_2$. One has a natural
$p$-prolongation sequence
$$(S^r_p[[{\bf T},\d_p {\bf T},...,\d_p^r{\bf T}]])_{r \geq 0}$$
(where $\d_p{\bf T},\d_p^2{\bf T},...$ are pairs of new variables). Then the $r+1$-tuple
$${\mathcal F},\d_p {\mathcal F},...,\d_p^r {\mathcal F}$$
defines a commutative formal group in $r+1$ variables $T,\d_p T,...,\d_p^r T$. Setting ${\bf T}=0$ in the above series, and forgetting about the first of them, we obtain an $r$-tuple of series
$$F_1:=\{\d_p {\mathcal F}\}_{|{\bf T}=0},...,F_r:=\{\d_p^r {\mathcal F}\}_{|{\bf T}=0}.$$
This $r$-tuple belongs to $S^r_p[\d_p {\bf T},...,\d_p^r {\bf T}]\h$ and defines a group
\begin{equation}
\label{groo}
(\hat{\mathbb A}^r_{S^r_p},[+])\end{equation}
  in the category of $p$-adic formal schemes over $S^r_p$.
Now let
$$l(T)=\sum_{n \geq 1} a_nT^n \in (S\otimes {\mathbb Q})[[T]]$$
be the logarithm of ${\mathcal F}$. Recall that $na_n \in S$ for all $n$. Define
\begin{equation}
\label{Lpr}
L^r_p:=\frac{1}{p}\{\phi^r(l(T))\}_{|T=0}\in (S^r_p\otimes {\mathbb Q})[[\d_p T,...,\d_p^r T]].\end{equation}
Then $L^r_p$ actually belong to $S^r_p[\d_p T,...,\d_p^r T]\h$ and define group homomorphisms
$$L^r_p:(\hat{\mathbb A}^r_{S^r_p},[+])\ra (\hat{\mathbb A}^1_{S^r_p},+)=\hat{\mathbb G}_{a,S^r_p}.$$
For all the facts above we refer to \cite{book}, pp. 123-125.

Now let $S^r_{\pi}:=\cO(J^r_{\pi}(X_{R_{\pi}})) \subset S^r_p \otimes_{R_p} R_{\pi}$; cf. Proposition \ref{inj}.
We have the following

\begin{proposition}
\label{xxx}
\

1) For some integer $n(r)\geq 1$, $p^{n(r)} F_r$  belongs to the image of the natural homomorphism
$$S^r_{\pi}[\d_{\pi} {\bf T},...,\d_{\pi}^r {\bf T}]\h \ra (S^r_{p}\otimes_R R[\pi])[\d_p {\bf T},...,\d_p^r {\bf T}]\h.$$

2) If $v_p(\pi)\geq \frac{1}{p-1}$, e.g. if $\pi=1-\zeta_p$, then $L_{\pi}^r:=\frac{p}{\pi}L^r_{p}$ belongs to the image of the natural homomorphism
$$S^r_{\pi}[\d_{\pi} T,...,\d_{\pi}^r T]\h \ra (S^r_{p}\otimes_R R[\pi])[\d_p T,...,\d_p^r T]\h.$$
\end{proposition}

{\it Proof}.
Since $\phi^r(T) \equiv T^{p^r}$ mod $\pi$ in
$R_{\pi}[\d_{\pi} T,...,\d_{\pi}^r T]$
we have $\{\phi^r(T)\}_{|T=0}\equiv 0$ mod $\pi$ in the same ring.
Set $G_{r,\pi}=\frac{1}{\pi} \{\phi^r(T)\}_{|T=0}$. We claim that for any $F \in S[[T]]$ with $F(0)=0$
we have
$$p^N \{\d_{p}^r F\}_{|T=0} \in S_{\pi}^r[\d_{\pi} T,...,\d_{\pi}^r T]\h$$
for some $N$.
Indeed since some power of $p$ times $\d_p^r F$ is a polynomial with ${\mathbb Z}$-coefficients
in $F,\phi(F),...,\phi^r(F)$ it is enough to show that $\{\phi^r(F)\}_{|T=0}$ is a restricted power series in $\d_{\pi} T,...,\d_{\pi}^r T$ for any $r$. But $\{\phi^r(F)\}_{|T=0}$ is a power series with $S^r_{\pi}$-coefficients in  $\{\phi^i(T)\}_{|T=0}=\pi G_{i,\pi}$, $i \leq r$, and our claim is proved. The same argument works for $T$ replaced by and tuple of variables; this ends the proof of assertion 1. To check assertion 2 note that
$$L_{\pi}^r=\sum_{n \geq 1} \phi(na_n) \frac{\pi^{n-1}}{n} G_{r,\pi}^n$$
and we are done because $na_n \in S$ and  $v_p(\pi^{n-1}/n)$ is $\geq 0$ and $\ra \infty$.
\qed

\subsection{Conjugate operators}
Recall from \cite{book}, Proposition 3.45, that if $X/R_p$ is a smooth affine
scheme and $\partial:\cO(X) \ra \cO(X)$
is an $R_p$-derivation then there are   unique
$R_p$-derivations
\begin{equation}
\label{cirip}
\partial_0,...,\partial_r:\cO(J^r_p(X))\ra \cO(J^r_p(X)),\end{equation}
called the {\it conjugate operators} of $\partial$ such that for all $s,j=0,...,r$:

i) $\partial_j \circ \phi^s=0$ on $\cO(X)$ if $j \neq s$;

2) $\partial_j \circ \phi^j=p^j \cdot \phi^j \circ \partial$ on $\cO(X)$.

\noindent In case $r=1$ it is trivial to see that

3) $\partial_1 f=0$ and $\partial_1  \d_p f=\phi \partial f$ for  $f \in \cO(X)$.

4) $\partial_0 f=\partial f$ and $\partial_0 \d_p f =-f^{p-1} \partial f$ for $f \in \cO(X)$.

 (By uniqueness this construction extends, in its obvious sheafified version, to the case when $X$ is not necessarily affine.) Note that the operator $\partial_r$  was introduced in \cite{Barcau} in a special case.
Clearly  $\partial_j$ uniquely extend to $R_{\pi}$-derivations $\partial_j:\cO(J^r_p(X))\otimes_{R_p} R_{\pi}\ra \cO(J^r_p(X))\otimes_{R_p} R_{\pi}$. View $\cO(J^r_{\pi}(X_{R_{\pi}}))$ as a subring of $\cO(J^r_p(X))\otimes_{R_p} R_{\pi}$. Then we have:

\begin{proposition}
\label{conjugate}
 If $f \in \cO(J_p^r(X))$ is $\d_{\pi}$-overconvergent then so are
 $\partial_j f$ for $j=0,...,r$.
Moreover if $r=1$ then $\cO(J^1_{\pi}(X_{R_{\pi}}))$ is sent into itself by both $\partial_0$ and $\partial_1$.
\end{proposition}

{\it Proof}. We may assume $X$ is affine. Write
$$f=F(f_1,...,f_n,...,\d_{\pi}^r f_1,...,\d_{\pi}^r f_n),\ \ F \in R_{\pi}[T_{10},...,T_{n0},...,T_{1r},...,T_{nr}]\h.$$
Then
$$\partial_j f=\sum_{is} \frac{\partial F}{\partial T_{is}}
(f_1,...,f_n,...,\d_{\pi}^r f_1,...,\d_{\pi}^r f_n)\partial_j (\d_{\pi}^s f_i).$$
So it is enough to show that  $\partial_j (\d_{\pi}^s f_i)$ times a power $p^{\nu}$ of $p$ is in $\cO(J^r_{\pi}(X_{R_{\pi}}))$; we also need to show that one can take $\nu=0$ if $r=1$. Since $\d_{\pi}^s f_i$ is a polynomial with $R_{\pi}$-coefficients in $\d_p f_i,...,\d_p^s f_i$ it is enough to show that $\partial_j \d_p^s f$ times a power $p^{\nu}$ of $p$ is in $\cO(J^r_{\pi}(X_{R_{\pi}}))$ and that one can take $\nu=0$ if $r=1$. These statements are clear from
the properties 1, 2, 3, 4 of the  the $\partial_j$s.
\qed

\section{Review of modular curves}

This section is entirely ``non-differential'' and represents a review of essentially well-known facts about modular parametrization and bad reduction of modular curves. These facts will play a role later in the paper.

\subsection{Modular parameterization}
Consider the following classes of objects:

\medskip

(1) Normalized newforms
\begin{equation}
\label{modf}
f=\sum_{n \geq 1} a_n q^n\end{equation}
 of weight $2$ on $\Gamma_0(M)$ over ${\mathbb Q}$; in particular $a_1=1$, $a_n \in {\mathbb Z}$.

\medskip

(2) Elliptic curves $A$ over ${\mathbb Q}$ of conductor $M$.

\medskip

Say that $f$ in (1) and $A$ in (2) correspond to each other if there exists
 a morphism
\begin{equation}
\label{shimura}
\Phi:X_0(M)\ra A
 \end{equation}
 over ${\mathbb Q}$  such that the pull back to $X_0(M)$ of some $1$-form on $A$ over ${\mathbb Q}$
corresponds to $f$ and $L(A,s)=\sum a_n n^{-s}$. We have the following fundamental result:

\begin{theorem}
\label{fundamental}\

 i) For any  $f$
as in (1) there exists an $A$ as in (2) which corresponds to $f$.

ii)
For any  $A$ as in (2) there exists $f$ as in (1) which corresponds to $A$.
\end{theorem}

The first part of the theorem is due to Eichler, Shimura, and Carayol; cf.  \cite{Knapp} for an exposition of this theory and  references. The second part of the Theorem is the content of the
 Tanyiama-Shimura conjecture proved, in its final form, in
\cite{four}.

\medskip

{\it From now on we fix $M=Np$ with $(N,p)=1$, $p \geq 5$, $N \geq 5$, and we fix $f$ and $A$  corresponding to each other, as in Theorem \ref{fundamental}. We further assume
 $a_p=1$ or, equivalently,  $A$ has split multiplicative reduction at $p$.}

\medskip

Recall from \cite{Knapp}, p.282, that
  $a_{mn}=a_m a_n$ for $(m,n)=1$, $a_{\ell^r} a_{\ell}= a_{\ell^{r+1}}+\ell a_{\ell^{r-1}}$ for $\ell \not| Np$,  and $a_{\ell^r}=a_{\ell}^r$ for $\ell | Np$; in particular $a_{p^r}=1$ for all $r$.

\subsection{Model of $X_1(Np)$ over ${\mathbb Z}[1/N,\zeta_p]$}
Recall  that the modular curve  $X_1(Np)$ over ${\mathbb C}$ has a model (still denoted by $X_1(Np)$ in what follows) over ${\mathbb Z}[1/N,\zeta_p]$
considered in  \cite{Gross}, p. 470;
this is a version of a curve  introduced by Deligne and Rappoport \cite{DR} and the two curves
become canonically  isomorphic over ${\mathbb Z}[1/N, \zeta_N, \zeta_p]$ if $\zeta_N$
is a fixed primitive $N$-th root of unity. Recall some of the main properties of $X_1(Np)$.
First $X_1(Np)$ is a regular scheme proper and flat  of relative dimension $1$ over ${\mathbb Z}[1/N,\zeta_p]$ and smooth over ${\mathbb Z}[1/Np,\zeta_p]$. Also the special fiber of $X_1(Np)$ over ${\mathbb F}_p$ is a union of two smooth projective curves $I$ and $I'$ crossing transversally at a finite set $\Sigma$ of points. Furthermore $I$ is isomorphic to the Igusa curve $I_1(N)$ in \cite{Gross}, p. 160, so $I$ is the smooth compactification of the curve classifying triples $(E,\alpha, \beta)$ with $E$ an elliptic curve over a scheme of characteristic $p$, and $\alpha:\mu_N \ra E$, $\beta:\mu_p \ra E$ are embeddings (of group schemes). Similarly  $I'$ is the smooth compactification of the curve classifying triples $(E,\alpha, b)$ with $E$ an elliptic curve over a scheme of characteristic $p$, and $\alpha:\mu_N \ra E$, $b:{\mathbb Z}/p{\mathbb Z}\ra E$  are embeddings.
Finally $\Sigma$ corresponds to the supersingular locus on the corresponding curves.

\subsection{N\'{e}ron model of $A$ over $R_{\pi}$}
Let $\pi=1-\zeta_p$ and consider a fixed embedding of ${\mathbb Z}[\zeta_N, \zeta_p, 1/N]$ into $R_{\pi}$ (hence of ${\mathbb Z}[\zeta_N,1/N]$ into $R_p$.]

Let $A_{R_{\pi}}$ be the Neron model of $A_{K_{\pi}}:=A\otimes_{\mathbb Q} K_{\pi}$ over $R_{\pi}$; cf. \cite{sil2}, p. 319. Then the $\pi$-adic completion
$(A_{R_{\pi}}^0)\h$
of the connected component $A_{R_{\pi}}^0$ of $A_{R_{\pi}}$ is isomorphic to the $\pi$-adic completion $({\mathbb G}_m)\h$ of ${\mathbb G}_m=Spec\ R_{\pi}[x,x^{-1}]$.
Indeed by \cite{sil2}, Theorem 5.3, p. 441,
$A_{K_{\pi}}$ is isomorphic over $K_{\pi}$ to a Tate curve $E_q/K_{\pi}$ with $q \in \pi R_{\pi}$. By \cite{sil2}, Corollary 9.1, p. 362,
$A_{R_{\pi}}^0$ is the smooth locus over $R_{\pi}$ of a projective curve defined by the
minimal Weierstrass equation of $A_{K_{\pi}}$. Now the defining Weierstrass
equation of the Tate curve (\cite{sil2}, p. 423) is already minimal (cf. \cite{sil2}, Remark 9.4.1, p. 364). The isomorphism $(A_{R_{\pi}}^0)\h\simeq ({\mathbb G}_m)\h$ then
follows from the formulae of the Tate parameterization \cite{sil2}, p. 425.

On the other hand recall that the modular curve $X_1(N)$  over ${\mathbb C}$ has a natural smooth projective model (still denoted by $X_1(N)$)
over ${\mathbb Z}[1/N]$ such that
$$Y_1(N):=X_1(N)\backslash (cusps)$$
 parameterizes pairs $(E,\alpha)$
consisting of elliptic curves $E$ with an embedding $\alpha:\mu_N \ra E$. The morphism $X_1(Np)\ra X_1(N)$ over ${\mathbb C}$  induces a
morphism
$$\epsilon: X_1(Np)_{R_{\pi}}\backslash \Sigma \ra X_1(N)_{R_{\pi}}\backslash (ss)$$
 over $R_{\pi}$, where $(ss)$ is the supersingular locus in the closed fiber of $X_1(N)_{R_{\pi}}$. Indeed the morphism $X_1(Np)\ra X_1(N)\ra J_1(N)$
over ${\mathbb C}$ (where $J_1(N)$ is the Jacobian of $X_1(N)$ over ${\mathbb C}$ and $X_1(N) \ra J_1(N)$ is the Abel-Jacobi map defined by $\infty$) induces a morphism
from $X_1(Np)_{R_{\pi}}\backslash \Sigma$ into the Jacobian $J_1(N)_{R_{\pi}}$ of  $X_1(N)_{R_{\pi}}$
(by the Neron property, because the latter Jacobian is an abelian scheme and hence is the Neron model of its generic fiber).
But the image of $X_1(Np)_{R_{\pi}}\backslash \Sigma \ra J_1(N)_{R_{\pi}}$ is clearly contained in the image of the Abel-Jacobi map
$X_1(N)_{R_{\pi}}\ra J_1(N)_{R_{\pi}}$ which gives a morphism $X_1(Np)_{R_{\pi}}\backslash \Sigma \ra X_1(N)_{R_{\pi}}$ and hence the desired morphism
$\epsilon:X_1(Np)_{R_{\pi}}\backslash \Sigma \ra X_1(N)_{R_{\pi}}\backslash (ss)$. Let $X \subset X_1(N)_{R_p} \backslash (ss)$ be an affine open set,
 $X_{R_{\pi}}:=X \otimes_{R_p} R_{\pi}\subset  X_1(N)_{R_{\pi}} \backslash (ss)$ its base change to $R_{\pi}$,
 and $X_!:=\epsilon^{-1}(X_{R_{\pi}})$. Denote by ${\mathcal X}_{R_{\pi}}$  the $\pi$-adic completion of $X_{R_{\pi}}$.
Also note that the $\pi$-adic completion of $X_!$ has two connected components; let ${\mathcal X}_!$ be the component whose reduction mod $\pi$ is contained in $I \backslash \Sigma$. We get a  morphism $\epsilon:{\mathcal X}_! \ra {\mathcal X}_{R_{\pi}}$.

\subsection{Igusa curve and lift to characteristic zero}
It will be useful to recall one of the possible constructions of the Igusa curve $I$. Let $L$ be the line bundle on $X_1(N)_{R_p}$
such that the sections of the powers of $L$ identify with the modular forms of various weights on $\Gamma_1(N)$; cf. \cite{Gross} p. 450 where $L$ was denoted by $\omega$. Let $E_{p-1}\in H^0(X_1(N)_{R_p},L^{p-1})$ be the normalized Eisenstein form of weight $p-1$ and let $(ss)$ be the supersingular locus on $X_1(N)_{R_p}$ (i.e. the zero locus of $E_{p-1}$).
(Recall that $E_{p-1}$ is normalized by the condition that its Fourier expansion has constant term $1$.)
Take an open  covering $(X_i)$ of $X$  such that  $L$ is trivial on each $X_i$ and we let $x_i$ be a basis of $L$ on $X_i$. Then $E_{p-1}=\varphi_i x_i^{p-1}$
where $\varphi_i \in \cO(X_i)$. Set $x_i=u_{ij}x_j$, $u_{ij}\in \cO^{\times}(X_{ij})$, $X_{ij}=X_i \cap X_j$. Consider the $R_{\pi}$-scheme
$X_{!!}$ obtained by gluing the schemes $X_{!!i}:=Spec\ \cO(X_{i,R_{\pi}})[t_i]/(t_i^{p-1}-\varphi_i)$ via $t_i=u_{ij}^{-1}t_j$ (where $X_{i,R_{\pi}}:=X_i\otimes_{R_p} R_{\pi}$). Note that $t_i^{p-1}-\varphi_i$ are monic polynomials whose derivatives are invertible in $\cO(X_{i,R_{\pi}})[t_i]/(t_i^{p-1}-\varphi_i)$.
Denote in the discussion below by an upper bar the functor $\otimes k$.
 Note that the scheme $\overline{X}_{!!}=X_{!!} \otimes k$ is isomorphic to $\overline{\mathcal X}_!={\mathcal X}_!\otimes k$; indeed
$\overline{X}_{!!}$ is clearly birationally equivalent to $I$ (cf. \cite{Gross}, pp. 460, 461) and is the integral closure of $\overline{X}$ in the fraction field
of $\overline{X}_{!!}$. We claim that:

 \begin{proposition}
 \label{liftt}
 The isomorphism $\overline{X}_{!!} \simeq \overline{\mathcal X}_!$  lifts uniquely to an isomorphism $(X_{!!})\h\simeq {\mathcal X}_!$.
 \end{proposition}

 {\it Proof}. Indeed this follows immediately by applying
 the  standard Lemma \ref{hensel} below to $S:=\cO({\mathcal X}_i)$, ${\mathcal X}_i=\widehat{X}_i$,
 $S_!=\cO({\mathcal X}_{!i})$, ${\mathcal X}_{!i}=\epsilon^{-1}({\mathcal X}_i)$.
 \qed

\begin{lemma}
 \label{hensel}
 Let $S \ra S_{!}$ be a morphism of flat $\pi$-adically complete $R_{\pi}$-algebras, let $f\in S[t]$ be a monic polynomial and assume we have an isomorphism $\overline{\sigma}:\overline{S}[t]/(\overline{f}) \ra \overline{S}_{!}$ such that $df/dt$ is invertible in $\overline{S}[t]/(\overline{f})$. Then $\overline{\sigma}$ lifts uniquely to an isomorphism
$\sigma:S[t]/(f)\ra S_{!}$.
\end{lemma}

{\it Proof}. The homomorphism $\sigma$ exists and is unique by Hensel's Lemma; it is an isomorphism because $\overline{\sigma}$ is one and $\pi$
is a non-zero divisor in both $S$ and $S_{!}$.
 \qed

 \subsection{Review of diamond operators}
 Recall from \cite{Gross} that $G:=({\mathbb Z}/p{\mathbb Z})^{\times}$ acts on the covering $X_1(Np) \ra X_1(N)$ over ${\mathbb Z}[1/N,\zeta_p]$
 via the diamond operators $\langle d \rangle_p$, $d \in G$; this action preserves the Igusa curve $I$ and induces on $I$ the usual diamond operators.
 In particular $I/G\ra X_1(N)_{{\mathbb F}_p}$ is an isomorphism. So $G$ acts on the covering $\epsilon:{\mathcal X}_!\ra {\mathcal X}_{R_{\pi}}$ and hence on the isomorphic covering $({\mathcal X}_{!!})\h\ra {\mathcal X}_{R_{\pi}}$; cf. Proposition \ref{liftt}. It is easy to explicitly find the latter action. Indeed
 any $G$-action on a covering $(X_{!!i})\h \ra X_{i,R_{\pi}}$ must have the form
 \begin{equation}
 \label{forrr}
 d \cdot t_i=\zeta_{p-1}^{\chi(d)}t_i,\ \ d \in G,\end{equation}
  for some homomorphism $\chi:G \ra {\mathbb Z}/(p-1){\mathbb Z}$, where $\zeta_{p-1}$ is a primitive root of unity of order $p-1$. Now we claim that $\chi$ must be an isomorphism. Indeed if $\chi$ was not surjective then
 the $G$-action on the Igusa curve $I$ would be such that $I/G\ra X_1(N)_{{\mathbb F}_p}$ has degree $>1$, a contradiction.

 \subsection{Classical  and $p$-adic modular forms}
We end by reviewing some more terminology and facts, to be used later, about classical modular forms and their relation with the $p$-adic modular forms of Serre and Katz.
 Let $M$ be any positive integer. (In applications we write $M=Np^{\nu}$, $(N,p)=1$.)
In what follows a {\it classical modular form} over a ring $B$, of weight $\kappa$, on $\Gamma_1(M)$
 will be understood in the sense of \cite{DR,Katz,Gross} as a rule that attaches to any $B$-algebra $C$
and any triple consisting of an elliptic curve  $E/C$, an embedding $\mu_{M,C} \ra
E[M]$, and an invertible one form on $E$ an element of $C$ satisfying the usual compatibility rules and the usual holomorphy condition for the Fourier expansion (evaluation on the Tate curve).
We denote by
$$M(B,\kappa,M)=M(B,\kappa,\Gamma_1(M))$$
the $B$-module of all these forms. We denote by
$$M(B,\kappa,\Gamma_0(M))$$
the submodule of those forms which are invariant under the usual
diamond operators.
In particular any newform as in (\ref{modf}) is an element of $M({\mathbb Z},2,\Gamma_0(Np))$; cf. \cite{DI}, p.113. Also by \cite{Gouvea}, p. 21,
the spaces $M(R_p,\kappa,Np^{\nu})$ embed into Katz's ring of {generalized $p$-adic modular forms} ${\mathbb W}={\mathbb W}(N,R_p)$ parameterizing trivialized elliptic curves $E$ over $p$-adically complete $R_p$-algebras, equipped with an embedding $\mu_N \subset E[N]$; if $f\in M(R_p,\kappa,Np^{\nu})$ then $f$, as an element of ${\mathbb W}$, satisfies $\lambda \cdot f=\lambda^{\kappa} f$ for $\lambda \in {\mathbb Z}_p^{\times}$, $\lambda \equiv 1$ mod $p^{\nu}$.
(Here $\lambda \cdot f$ denotes the action of ${\mathbb Z}_p^{\times}$ on ${\mathbb W}$ induced by changing the trivialization.) If $f$ is actually in $M(R_p,\kappa,\Gamma_0(Np^{\nu}))$ then $\lambda \cdot f=\lambda^{\kappa} f$ for $\lambda \in {\mathbb Z}_p^{\times}$. In particular any newform $f$ as in (\ref{modf})
on $\Gamma_0(Np)$ defines an element (still denoted by $f$) of ${\mathbb W}$ such that $\lambda \cdot f=\lambda^2f$, $\lambda\in {\mathbb Z}^{\times}_p$. By \cite{Gouvea}, p.21, any  $f$ as in (\ref{modf}) on $\Gamma_0(Np)$ is a $p$-adic modular form of weight $2$ over $R_p$ {\it in the sense of Serre}, i.e. it is a $p$-adic limit in ${\mathbb W}$ (or equivalently in $R_p[[q]]$) of classical modular forms over $R_p$ of weight $\kappa_n\in {\mathbb Z}$ on $\Gamma_1(N)$ with $\kappa_n \equiv 2$ mod $p^n(p-1)$.
So if $f=\sum a_n q^n$ is as in (\ref{modf}) on $\Gamma_0(Np)$ then, by \cite{Gouvea}, p. 36, $\sum a_n q^{np}$ is also a $p$-adic modular form of weight $2$ in the sense of Serre. In particular the reduction mod $p$ of
$\sum a_n q^{np}$ is the expansion of a modular form over $k$ on $\Gamma_1(N)$ of weight $\equiv 2$ mod $p-1$. Finally recall from \cite{Gross} that the Serre operator $\theta:=d\frac{d}{dq}:k[[q]]\ra k[[q]]$ increases weights of classical modular forms over $k$ by $p+1$.
We conclude that the image in $k[[q]]$  of
$$\sum_{(n,p)=1} \frac{a_n}{n} q^n\in R_p[[q]]$$
 is the expansion of a modular form over $k$ on $\Gamma_1(N)$ of weight $\equiv 0$ mod $p-1$.

We end by recalling a few basic facts about Hecke operators.
Throughout the discussion below the divisors  of a given non-zero integer  are always taken to be positive, the greatest common divisor of two non-zero integers $m,n$ is denoted by $(m,n)$, and we use the convention $(m,n)=n$ for $m=0$, $n \neq 0$. Fix again  a positive integer  $M$  and let $\epsilon_{M}:{\mathbb Z}_{> 0}\ra \{0,1\}$ be the ``trivial primitive character" mod $M$ defined by $\epsilon_{M}(A)=1$ if $(A,M)=1$ and $\epsilon_{M}(A)=0$ otherwise.

For each integers $n \geq 1$, $\k \geq 2$ and any ring $C$
  define the operator  $T_{\k, M}(n):C[[q]] \ra C[[q]]$
by the formula
$$
T_{\k,M}(n)f= \sum_{m \geq 0}\left(
\sum_{A|(n,m)} \epsilon_{M}(A)A^{\k-1}a_{\frac{mn}{A^2}}\right) q^m.$$
Recall (cf., say, \cite{Knapp})
that if $f=\sum_{m \geq 0}a_mq^m \in {\mathbb C}[[q]]$ is the Fourier expansion of a form in $M({\mathbb C},\kappa,\Gamma_0(M))$,  $\k\geq 2$, then
the series $T_{\k,M}(n)f$ is the Fourier expansion of the corresponding Hecke operator on $f$.
 Note that
if $M=Np^{\nu}$, $(N,p)=1$, $(n,p)=1$ then  $T_{\k,N}=T_{\k,M}$ as operators on $C[[q]]$. Now if $f$ is as in (\ref{modf}) then $T_{2,Np}(n)f=a_n f$ for all $n \geq 1$; so, for $(n,p)=1$ we have $T_{2,N}(n)f=a_nf$. On the other hand, going back to an arbitrary $f =\sum a_m q^m \in C[[q]]$, we have
$$T_{\k,N}(p)f=\sum_m a_{mp}q^m + p^{\kappa -1}\sum_m a_m q^{pm},$$
$$T_{\k,Np}(p)f=\sum_m a_{mp}q^m.$$
So $T_{\k,N}(p)\equiv T_{\k,Np}(p)$ mod $p$ as operators on $C[[q]]$. Specializing again to $f\in {\mathbb Z}[[q]]$ as in (\ref{modf}) on $\Gamma_0(Np)$ we have
$T_{2,Np}(p)f=a_p f=f$
so we get $T_{2,N}(p)f\equiv f$ mod $p$ in ${\mathbb Z}[[q]]$.

\section{$\d_{p}$-modular forms arising from bad reduction}

In this section we return to ``differential matters''. We will use bad reduction of the modular curve $X_1(Np)$ at $p$ to construct certain $\d_{\pi}$-functions
on this curve.
These functions will then induce (via a {\it geometric trace} construction)  certain new interesting $\d_{\pi}$-{\it modular forms}
 on the modular curve $X_1(N)$. By further applying the {\it arithmetic trace}
from $R_{\pi}$ down to $R_p$ we will obtain certain new $\d_p$-modular forms
 on $X_1(N)$.
We will then analyze the $\d_{\pi}$-{\it Fourier expansions} (respectively $\d_p$-Fourier expansions) of these forms. On our way of doing this we will review the concepts of $\d_p$-{\it modular form} and $\d_p$-{\it Fourier expansion} following \cite{difmod,book}.

\subsection{$\d_p$-modular forms and $\d_{\pi}$-modular forms}
Let $L$ be the line bundle on $X_1(N)_{R_p}$
such that the spaces of sections $H^0(X_1(N)_{R_p},L^{\otimes \kappa})$ identify with the spaces
$M(R_p,\kappa,N)$ of classical modular forms over $R_p$ of  weight $\kappa$ on $\Gamma_1(N)$; cf. \cite{Gross} p. 450 where $L$ was denoted by $\omega$.

Let $X\subset X_1(N)_R$ an affine open subset.
(In \cite{book,igusa, hecke} we always assumed that $X$ is disjoint from the cusps; we will not assume this here because we find it convenient  to cover a slightly more general case.)
The restriction of $L$ to $X$ will still be denoted by $L$. Consider the $X$-scheme
\begin{equation}
\label{defV}
V:=Spec\left( \bigoplus_{n \in \bZ} L^{\otimes n}\right).\end{equation}
By a {\it $\d_{p}$-modular function of order $r$ on $X$} \cite{igusa} we understand an element of the ring $M^r_p:=\cO(J^r_p(V))$.
If we set, as usual,  $V_{R_{\pi}}:=V \otimes_{R_p} R_{\pi}$ then by a {\it $\d_{\pi}$-modular function of order $r$ on $X$} we will understand an element of
$M^r_{\pi}:=\cO(J^r_{\pi}(V_{R_{\pi}}))$.
 The formation of these rings is functorial in $X$. Also if $L$ is trivial on $X$ with basis $x$ then $M^r_p$ identifies with
 $\cO(J^r_p(X))[x,x^{-1},\d_p x,...,\d_p^r x]\h$ and $M^r_{\pi}$ identifies with  $\cO(J^r_{\pi}(X))[x,x^{-1},\d_{\pi} x,...,\d_{\pi}^r x]\h$.
Recall the ring $W:={\mathbb Z}[\phi]$ of polynomials in $\phi$; it will play in what follows the role of ring of {\it weights}.
 By a {\it $\d_p$-modular form of order $r$ and  weight $w\in W$ on $X$} we mean a $\d_p$-modular function $f \in M^r_p$ such that  for each $i$, $f\in \cO(J^r_p(X_i))\cdot x_i^w$; cf. \cite{igusa}.
We denote by $M^r_p(w)$ the $R_p$-module of $\d_p$-modular forms of order $r$ and  weight $w$ on $X$.
For $w=0$ we set $S_p^r=M^r_p(0)=\cO(J^r_p(X))$.
By a {\it $\d_{\pi}$-modular form of order $r$ and weight $w$ on $X$} we will mean a $\d_{\pi}$-modular function $f \in M^r_{\pi}$ such that  for each $i$, $f\in \cO(J^r_{\pi}(X_{i,R_{\pi}}))\cdot x_i^w$. We denote by $M^r_{\pi}(w)$ the $R_{\pi}$-module of $\d_{\pi}$-modular forms of order $r$ and  weight $w$ on $X$. For $w=0$ we set $S_{\pi}^r=M^r_{\pi}(0)=\cO(J^r_{\pi}(X_{R_{\pi}}))$. In view of (\ref{tra}) and (\ref{tr}) we have natural $R_{\pi}$-algebra  homomorphisms
\begin{equation}
\label{uncal}
M^r_{\pi}\ra M^r_p \otimes_{R_p} R_{\pi}\end{equation}
preserving weights, i.e. inducing $R_{\pi}$-linear maps
$$M^r_{\pi}(w)\ra M^r_p(w)\otimes_{R_p} R_{\pi},\ \ w\in W.$$
Also we have  $R_p$-linear trace maps
\begin{equation}
\label{doicai}
\tau_{\pi}:M^r_{\pi} \ra
M^r_p\end{equation}
 that preserve weights i.e. induce maps
\begin{equation}
\label{vopsea}
\tau_{\pi}:M^r_{\pi}(w) \ra M^r_p(w), \ \ w \in W.\end{equation}
In particular we have $R_{\pi}$-algebra  homomorphisms
$$S^r_{\pi}\ra S^r_p\otimes_{R_p} R_{\pi}$$ and $R$-linear trace maps
$$\tau_{\pi}:S^r_{\pi}\ra S_p^r.$$

When applied to the scheme $V$, Definition \ref{basicdef} translates into the following:

\begin{definition}
A $\d_p$-modular function $f \in M^r_p$ is called $\d_{\pi}$-{\it overconvergent}
if one of the following equivalent conditions is satisfied:

1) $f$ times a power of $p$ belongs to the image of the map (\ref{uncal});

 2) $f$ times a power of $p$  belongs  to the image of the map (\ref{doicai}).
\end{definition}

 \subsection{$\d_{\pi}$-modular forms from $\d_{\pi}$-functions on ${\mathcal X}_!$}
Let $X \subset X_1(N)_{R_p}$ be disjoint from the supersingular locus $(ss)$
(but necessarily from $(cusps)$ !).
There is a canonical way of constructing $\d_{\pi}$-modular forms
of weights $0,-1,...,-p+2$ on $X$
from $\d_{\pi}$-functions on ${\mathcal X}_!$. Indeed we will construct
   natural {\it geometric trace} maps
 \begin{equation}
 \label{mappp}
 \tau_{\kappa}:\cO(J^r_{\pi}({\mathcal X}_!)) \ra  M^r_{\pi}(-\kappa),\ \ \ \kappa=0,...,p-2,\end{equation}
as follows. The isomorphism $(X_{!!})\h\simeq {\mathcal X}_!$ in Proposition \ref{liftt} induces an isomorphism $J^r_{\pi}({\mathcal X}_!)\simeq
J^r_{\pi}(X_{!!})$. Since $X_{!!i}:=Spec\ \cO(X_{i,R_{\pi}})[t_i]/(t_i^{p-1}-\varphi_i)$ is \'{e}tale over $X_{i,R_{\pi}}$ and since the formation of $\pi$-jet spaces commutes with \'{e}tale maps it follows that we have an identification
\begin{equation}
\label{iddee}
\cO(J^r_{\pi}(X_{!!i}))=\cO(J^r_{\pi}(X_{i,R_{\pi}}))[t_i]/(t_i^{p-1}-\varphi_i).\end{equation}
Let us denote the class of $t_i$ in the latter ring again by $t_i$ and let the image of $\alpha \in \cO(J^r_{\pi}({\mathcal X}_!))\simeq
\cO(J^r_{\pi}(X_{!!}))$ in $\cO(J^r_{\pi}(X_{!!i}))$ be $\sum_{\kappa=0}^{p-2} \alpha_{\kappa,i} t_i^{\kappa}$, $\alpha_{\kappa,i} \in \cO(J^r_{\pi}(X_{i,R_{\pi}}))$. Then define
$$\tau_{\kappa,i} \alpha:=\alpha_{\kappa,i} x_i^{-\kappa} \in \cO(J^r_{\pi}(X_{i,R_{\pi}})) \cdot x_i^{-\kappa}.$$
Note that from the equalities
$$\sum_{\kappa=0}^{p-2} \alpha_{\kappa,j} t_j^{\kappa}=\sum_{\kappa=0}^{p-2} \alpha_{\kappa,i} t_i^{\kappa}=
\sum_{\kappa=0}^{p-2} \alpha_{\kappa,i} u_{ij}^{-\kappa} t_j^{\kappa}$$
it follows that $\alpha_{\kappa,i}=u_{ij}^{\kappa}\alpha_{\kappa,j}$ hence
 $\tau_{\kappa,i} \alpha=\tau_{\kappa,j}\alpha$ for all $i$ and $j$. So the latter give rise to well defined  elements $\tau_{\kappa} \alpha \in M^r_{\pi}(-\kappa)$ which ends the construction of the map (\ref{mappp}).

 \begin{proposition}
 The map
 \begin{equation}
 \cO(J^r_{\pi}({\mathcal X}_!)) \ra  \bigoplus_{\kappa=0}^{p-2} M^r_{\pi}(-\kappa),\ \ \ \alpha \mapsto (\tau_0 \alpha,...,\tau_{p-2}\alpha)\end{equation} is an isomorphism.
 \end{proposition}

 {\it Proof}.
 Injectivity is clear from construction. Surjectivity immediately follows by reversing the construction
 of the trace maps above.
 \qed.

\medskip

On the other hand it will be useful to have a criterion
saying when a $\d_{\pi}$-function on ${\mathcal X}_!$ ``comes from'' a $\d_{\pi}$-modular form on $X$ of weight $0$, i.e. from a $\d_{\pi}$-function on $X$. Indeed recall the $G=({\mathbb Z}/p{\mathbb Z})^{\times}$-action on ${\mathcal X}_!$ induced by the diamond operators.
This action induces a $G$-action on $\cO(J^r_{\pi}({\mathcal X}_!))$ for all $r \geq 1$. Then we have:

\begin{proposition}
\label{invjjj}
The ring $\cO(J^r_{\pi}({\mathcal X}_!))^G$ of $G$-invariant elements of  $\cO(J^r_{\pi}({\mathcal X}_!))$ equals
$\cO(J^r_{\pi}(X_{R_{\pi}}))$.
\end{proposition}

{\it Proof}. This follows immediately from the identification (\ref{iddee}) and the fact that $G$ acts on $t_i$ by the formula
(\ref{forrr}) where $\chi$ is an isomorphism.
\qed

\subsection{$\d_{\pi}$-character composed with modular parameterization}
We assume, unless otherwise specified, that $\pi=1-\zeta_p$ and we fix, as usual an embedding ${\mathbb Z}[1/N,\zeta_N, \zeta_p]\ra R_{\pi}$. Also recall our fixed elliptic curve $A$ with modular parame\-triza\-tion (\ref{shimura}) and the modular form $f$ in (\ref{modf}). We continue to consider $X \subset X_1(N)_{R_p}$  an affine open set disjoint from $(ss)$.
We shall freely use the notation in our section on bad reduction.
By the N\'{e}ron property \cite{sil2}, p. 319, we get a morphism  $\Phi:X_! \ra A_{R_{\pi}}$ over $R_{\pi}$.
We get an induced morphism from ${\mathcal X}_!$
 into the connected component  $(A^0_{R_{\pi}})\h\simeq ({\mathbb G}_m)\h$.
This morphism $\Phi^0:{\mathcal X}_! \ra ({\mathbb G}_m)\h$ induces a morphism $\Phi^1:J^1_{\pi}({\mathcal X}_!)\ra J^1_{\pi}({\mathbb G}_m)$.
Now take the standard $\d_{\pi}$-character
$\psi_{\pi}\in \cO(J^1_{\pi}({\mathbb G}_m))$, cf. (\ref{sol}),
 identified with a morphism
$\psi_{\pi}:J^1_{\pi}({\mathbb G}_m)\ra ({\mathbb A}^1_{R_{\pi}})\h$.
By composition we get an induced morphism
 $f^{\sharp}:=\psi_{\pi}\circ \Phi^1:J^1_{\pi}({\mathcal X}_!) \ra ({\mathbb A}^1_{R_{\pi}})\h$.
 This morphism can be identified with an element
 \begin{equation}
 \label{defoffshhh}
 f^{\sharp}_{\pi}\in \cO(J^1_{\pi}({\mathcal X}_!)).\end{equation}
  (Here $f$ in $f_{\pi}^{\sharp}$ refers to the
 newform $f=\sum a_n q^n$ (\ref{modf}).)
Now, since $f$ is a form on $\Gamma_0(Np)$ it follows that $\Phi:X_!\ra A_{R_{\pi}}$ is invariant under the diamond operators
$\langle d \rangle_p$, $d \in G$. This implies that $f^{\sharp}_{\pi}$ is $G$-invariant. By Proposition \ref{invjjj} it follows that
\begin{equation}
 \label{defoffsh}
 f^{\sharp}_{\pi}\in \cO(J^1_{\pi}(X_{R_{\pi}}))=M^1_{\pi}(0)=S_{\pi}^1\end{equation}
i.e. $f^{\sharp}_{\pi}$ is a $\d_{\pi}$-modular form of weight $0$. Consequently its
image via the corresponding map  (\ref{vopsea}) defines a $\d_p$-modular form of weight $0$,
\begin{equation}
\label{simainou}
\tau_{\pi}  f_{\pi}^{\sharp} \in M^1_p(0)=S^1_p.\end{equation}

\subsection{$\d_p$-Fourier expansions and $\d_{\pi}$-Fourier expansions}

The $R_p$-point $\infty$ on $X_1(N)_{R_p}$ induces
$\d_{\pi}$-{\it Fourier expansion maps}
$$E:\cO(J^r_{\pi}({\mathcal X}_!))\ra R_{\pi}((q))[\d_{\pi} q,...,\d_{\pi}^r q]\h.$$
Indeed to construct such a map we may assume $X$ contains $\infty$; but in this case the map arises because $X_{!!}\ra X_{R_{\pi}}$ is \'{e}tale so the inverse image of $\infty$ by this map is a disjoint union of $R_{\pi}$-points.

On the other hand there are $\d_{\pi}$-{\it Fourier expansion maps}
\begin{equation}
\label{firstE}
E:M^r_{\pi}\ra R_{\pi}((q))[\d_{\pi} q,...,\d_{\pi}^r q]\h.\end{equation}
compatible, in the obvious sense, with the previous ones and with the $\d_p$-Fourier expansion maps in \cite{book, igusa}
\begin{equation}
\label{secondE}
E:M^r_p  \ra R_p((q))[\d_p q,...,\d_p^r q]\h.\end{equation}
We recall \cite{book} the $\d_p$-Fourier {\it expansion principle} according to which for any $w$ the map
$$E:M^r_p(w)\ra R_p((q))[\d_p q,...,\d_p^r q]\h$$
is injective and has a torsion free cokernel.

\begin{remark}
The maps (\ref{firstE}) and (\ref{secondE}) commute with the trace maps
$\tau_{\pi}: M^r_{\pi}\ra M^r_p$ and $\tau_{\pi}:R_{\pi}((q))[\d_{\pi} q,...,\d_{\pi}^r q]\h\ra
R_p((q))[\d_p q,...,\d_p^r q]\h$, in the sense that $E \circ \tau_{\pi}=\tau_{\pi} \circ E$.
\end{remark}

\begin{remark}
Clearly if $f \in M^r_p$ of $\d_{\pi}$-overconvergent then its $\d_p$-Fourier expansion $E(f)$ is $\d_{\pi}$-overconvergent. Later we will prove the $\d_{\pi}$-overconvergence of a number of remarkable $\d_p$-modular functions. By the present remark we will also get that their $\d_p$-Fourier expansions are $\d_{\pi}$-overconvergent. However the $\d_{\pi}$-overconvergence of all these expansions can also be proved directly.
\end{remark}

The next Proposition establishes a link between the $\d_{\pi}$-Fourier expansions of $\d_{\pi}$-functions on ${\mathcal X}_!$ and $\d_{\pi}$-Fourier
expansions of their geometric traces. Recall the series $E_{p-1}(q):=E(E_{p-1})\in R_p[[q]]$ and the fact that $E_{p-1}(q) \equiv 1$ mod $p$ in $R_p[[q]]$ \cite{Katz}.
So the series $E_{p-1}(q)$ has a unique $(p-1)$-root $\epsilon(q) \in R_p[[q]]$ such that $\epsilon(q) \equiv 1$ mod $p$ in $R_p[[q]]$.

\begin{proposition}
\label{forr}
If $\alpha \in \cO(J^1_{\pi}({\mathcal X}_!))$ then its $\d_{\pi}$-Fourier expansion is given by
$$E(\alpha)=\sum_{\kappa=0}^{p-2} E(\tau_{\kappa} \alpha) \epsilon(q)^{\kappa}.$$
\end{proposition}

{\it Proof}.
Shrinking $X$ we may assume $X=X_i$ for some $i$. From $E_{p-1}=\varphi_i x_i^{p-1}$ we get
$$E_{p-1}(q)=E(\varphi_i) E(x_i)^{p-1}=E(t_i)^{p-1} E(x_i)^{p-1}.$$
So $E(t_i x_i)=c \cdot \epsilon(q)$, $c \in R_p^{\times}$, $c^{p-1}=1$. Now the birational isomorphism between the Igusa curve $I$
and $\overline{X}_{!!}$ sends $t_ix_i$ into the form $a$ in \cite{Gross}, p. 460-461, and the Fourier expansion in $k[[q]]$ of the form $a$ at $\infty$ is $1$. It follows that $c=1$. We get
$$\begin{array}{rcl}
E(\alpha) & = & \sum_{\kappa=0}^{p-2} E(\alpha_{\kappa,i}) E(t_i)^{\kappa}\\
\  & \  & \  \\
\  & = & \sum_{\kappa=0}^{p-2} E(\alpha_{\kappa,i}) E(x_i)^{-\kappa} E(x_i)^{\kappa} E(t_i)^{\kappa}\\
\  & \  & \  \\
\  & = & \sum_{\kappa=0}^{p-2} E(\tau_{\kappa} \alpha) \epsilon(q)^{\kappa}.
\end{array}$$
\qed

\begin{proposition}
\label{doi}
If $f_{\pi}^{\sharp}\in \cO(J^1_{\pi}({\mathcal X}_!))$ is attached to $f=\sum a_n q^n$ an in (\ref{defoffsh}) then its $\d_{\pi}$-Fourier expansion
$E(f_{\pi}^{\sharp})\in R_{\pi}[[q]][\d_{\pi}q]\h$
has the form:
\begin{equation}
\label{mess}
\begin{array}{rcl}
E(f^{\sharp}_{\pi}) & = & \frac{1}{\pi}(\phi-p)\sum_{n \geq 1} \frac{a_n}{n} q^n\\
  & \  & \  \\
  \  & = &
\frac{1}{\pi} \left[ \left(\sum_{n \geq 1} \frac{a_n}{n}(q^p+\pi\d_{\pi}q)^n\right) -p \left( \sum_{n \geq 1} \frac{a_n}{n}q^n\right)\right]\\
\  & \  & \ \\
\  & = &
\frac{1}{\pi} \left[ \left(\sum_{n \geq 1} \frac{a_n}{n}(q^p+p\d_p q)^n\right) -p \left( \sum_{n \geq 1} \frac{a_n}{n}q^n\right)\right].
\end{array}\end{equation}
\end{proposition}

{\it Proof}. Entirely similar to the proof of Theorem 6.3 in \cite{eigen}
\qed

\begin{remark}
The series in the right hand side of Equation \ref{mess} are a priori elements of
$$K_{\pi}[[q,\d_{\pi}q]]=K_{\pi}[[q,\d_p q]].$$
The Lemma says in particular that these series are actually in  $R_{\pi}[[q]][\d_{\pi}q]\h$. One can also check the latter directly.
\end{remark}

\begin{proposition}
\label{majj}
The form $\tau_{\pi} f^{\sharp}_{\pi}$ in (\ref{simainou})
 satisfies the following
identity in the ring $R_p[[q]][\d_p q]\h$:
$$ E(\tau_{\pi}f^{\sharp}_{\pi}) =
\frac{p-1}{2}\left[ \left(\sum_{n \geq 1} \frac{a_n}{n}(q^p+p\d_p q)^n\right) -p \left( \sum_{n \geq 1} \frac{a_n}{n}q^n\right)\right].$$
\end{proposition}

{\it Proof}.
This follows from Proposition  \ref{doi} by using $Tr(\frac{1}{\pi})=\frac{p-1}{2}$.
\qed

\medskip

One can get a more explicit picture mod $\pi$ (respectively mod $p$) as follows.

\begin{proposition}
\label{trei}
The form $f_{\pi}^{\sharp}$ in  (\ref{defoffsh}) satisfies the following congruence mod $\pi$ in the ring $R_{\pi}[[q]][\d_{\pi}q]\h$:
$$ E(f_{\pi}^{\sharp})\equiv \left(\sum_{n\geq 1} a_n q^{np}\right) \cdot \frac{\d_{\pi}q}{q^p}-
\left(\sum_{n\geq 1} a_n q^{np^2}\right) \cdot \left(\frac{\d_{\pi}q}{q^p}\right)^p.$$
\end{proposition}

{\it Proof}.
Using  Proposition \ref{doi}
and the fact that $a_{mn}=a_ma_n$ for $(m,n)=1$ and $a_{p^i}=1$ for all $i$ \cite{Knapp}, p.282, one gets
immediately  that
$$E(f_{\pi}^{\sharp})_{|\d_{\pi}q=0}=-\frac{p}{\pi} \sum_{(m,p)=1} \frac{a_m}{m}q^m \equiv 0\ \ \ mod\ \ \pi.$$
Also the coefficient of the monomial $q^{p(n-1)}\d_{\pi}q$ in $E(f_{\pi}^{\sharp})$ equals $a_n$.
Finally fix  $i \geq 2$; the coefficient of the monomial $q^{p(n-i)}(\d_{\pi}q)^i$ in $E(f_{\pi}^{\sharp})$
equals
$$c_{i,n}:=\frac{\pi^{i-1}}{i!} (n-1)(n-2)...(n-i+1) a_n\in K_{\pi}.$$
If $i<p$ clearly $v_p(c_{i,n})>0$. If $i>p$ or if $i=p$ and $(n,p)=1$ then
$$v_p((n-1)(n-2)...(n-i+1))\geq 1$$
and since
$$v_p\left(\frac{\pi^{i-1}}{i!}\right)\geq\frac{i-1}{p-1}-\frac{i}{p-1}=-\frac{1}{p-1}.$$
we get $v_p(c_{i,n})>0$. Finally, assume $i=p$ and $p|n$. Then
$$c_{i,n}\equiv \frac{\pi^{p-1}}{p}a_n\equiv -a_n \ \ mod\ \ \pi$$
because
$$\begin{array}{rcl}
p & = & \pi^{p-1}(1+\zeta_p)(1+\zeta_p+\zeta_p^2)...(1+\zeta_p+...+\zeta_p^{p-2})\\
\  & \  & \  \\
\  & \equiv & \pi^{p-1}(p-1)!\ \ \ mod\ \ \ \pi\\
\  & \  & \  \\
\  & \equiv & -\pi^{p-1}\ \ \ mod\ \ \ \pi
\end{array}$$
which easily concludes the proof because $a_{np}=a_n$.
\qed

\begin{proposition}
\label{patru}
The form $\tau_{\pi} f_{\pi}^{\sharp}\in S^1_p$ in (\ref{simainou})
belongs to $pS^1_p$. Moreover the form $$f_p^{\sharp}:=\frac{2}{p}\tau_{\pi}f^{\sharp}_{\pi}\in S^1_p$$
is $\d_{\pi}$-overconvergent and
 satisfies the following
congruence mod $p$
$$ E(f_{p}^{\sharp})\equiv
  \left(\sum_{(n,p)=1}\frac{a_n}{n}q^n\right)-
\left(\sum_{n \geq 1} a_n q^{np}\right) \frac{\d_p q}{q^p}
$$
in the ring $R_p[[q]][\d_p q]\h$.
\end{proposition}

{\it Proof}.
By Proposition \ref{majj},
 one gets
$$E(\tau_{\pi}f^{\sharp}_{\pi})_{|\d_p q=0}=-\frac{p(p-1)}{2}\sum_{(n,p)=1}\frac{a_n}{n}q^n.$$
The coefficient of $q^{p(n-1)}\d_p q$ in $E(\tau_{\pi}f^{\sharp}_{\pi})$ equals
$pa_n$.
Also, for $i \geq 2$, the coefficient of $q^{p(n-i)}(\d_p q)^i$ in
$E(\tau_{\pi}f^{\sharp}_{\pi})$ equals
$$\frac{p-1}{2} \frac{p^i}{i!} a_n (n-1)(n-2)...(n-i+1).
$$
In particular $E(\tau_{\pi}f^{\sharp}_{\pi})$ is divisible by $p$ in the ring
$R[[q]][\d_p q]\h$. By the $\d_p$-Fourier expansion principle it follows that
$\tau_{\pi}f^{\sharp}_{\pi}$ is divisible by $p$ in $S^1_p$ which proves the first assertion of the Proposition. $\d_{\pi}$-overconvergence follows from
Proposition \ref{cconverse}. The rest of the Proposition
 then follows from the above coefficient computations.
\qed

\begin{remark}
 Let $\overline{f}=\sum a_m q^m\in k[[q]]$, $\overline{f^{(-1)}}:=\sum_{(n,p)=1}\frac{a_n}{n}q^n\in k[[q]]$ and let $V$ be $k$-algebra endomorphism of $k[[q]]$ that sends $q$ into $q^{p}$.
Then the series in $k[[q]]$ obtained from the right hand side of the formula in Proposition \ref{patru} by  reducing mod $p$ equals
$$\overline{g}:= \overline{f^{(-1)}} - V(\overline{f}) \frac{\d_p q}{q^p} \in k[[q]][\d_p q].$$
This series $\overline{g}$
is {\it Taylor $\d_p$-$p$-symmetric} in the sense of \cite{hecke}. Also, recalling from \cite{hecke} the  operators
denoted by $``pU"$ and $``pT_0(p)"$ acting on Taylor $\d_p-p$-symmetric series  and using
the fact that $T_{2,N} (p)\overline{f}=  \overline{f}$
it is a trivial exercise (using the formulae in \cite{hecke}) to check that
$``pU" \overline{g}=\overline{g}$ and hence
$$``pT_0(p)" \overline{g}=\overline{g}+V(\overline{f^{(-1)}}).$$
In particular note that $\overline{g}$ is {\it not} an eigenvector of $``pT_0(p)"$.
On the other hand
an action of the operators $T_0(n)$ (for level $N$)
on $k[[q]][\d_p q]$ was introduced in \cite{hecke};
using the fact that $T_{2,N}(n)\overline{f}=a_n \overline{f}$ for $(n,p)=1$ it follows (using the formulae in \cite{hecke}) that $nT_0(n)\overline{g}=a_n \overline{g}$ for $(n,p)=1$. So $\overline{g}$ {\it is} an eigenvector
of all operators $nT_0(n)$ with eigenvalues $a_n$.
\end{remark}

\section{$\d_{\pi}$-overconvergence of some basic $\d_p$-modular forms}

In this section we prove the $\d_{\pi}$-overconvergence of some of the basic $\d_p$-functions of the theory in \cite{char,difmod,Barcau,book,eigen,BP}.

\subsection{Review of the $\d_p$-modular forms $f^r_p$ \cite{book}}
We start by reviewing the construction of some basic $\d_p$-modular forms
$f^r_p=f^r_{p,jet}\in M^1_p(-1-\phi^r)$, $r \geq 1$. These were introduced in \cite{difmod, book}.
(There is a ``crystalline definition" of these forms introduced in \cite{difmod} for $r=1$ and \cite{Barcau} for $r \geq 1$ in the case of level $1$, and  in \cite{book} for arbitrary level; the equivalence of these definitions follows from  \cite{book}, Proposition 8.86.)
Below we follow \cite{book}, p. 263. The construction is as follows.
We let $X \subset X_1(N)_{R_p}$ be an affine open set disjoint from $(cusps)$.
Assume first that $L$ is trivial on $X$ and let $x$ be a basis of $L$.
Consider the universal elliptic curve $E \ra X$ over $R_p$ and view $x$ as a relative $1$-form on $E/X$. Cover $E$
by affine open sets $U_i$. Then the natural projections $J_p^r(U_i) \ra \hat{U}_i \hat{\otimes}_{S^0_p} S^r_p$ possess sections
$$s_{i,p}:\hat{U}_i \hat{\otimes}_{S^0_p} S^r_p \ra J_p^r(U_i).$$
Let $N^r_p:=Ker(J^r_p(E) \ra \hat{E}\otimes_{S_p^0} S_p^r)$; it is a group object in the category of $p$-adic formal schemes over $S_p^r$. Then the differences $s_{i,p}-s_{j,p}$ define
morphisms
$$s_{i,p}-s_{j,p}: \hat{U}_{ij} \hat{\otimes}_{S^0_p} S^r_p \ra N_p^r$$
where the difference is taken in the group law of $J_p^r(E)/S_p^r$.
On the other hand $N_p^r$ identifies with the group $(\hat{\mathbb A}^r_{S_p^r},[+])$ in
(\ref{groo})
with coordinates given by the $\d_p T,...,\d_p^r T$, where  $T$ is a parameter at the origin of $E$ chosen such that
$x \equiv dT$ mod $T$. Let $L_p^r$ be the series in (\ref{Lpr}) attached to the formal group of $E$
with respect to the same parameter $T$, viewed as a homomorphism $L_p^r:N^r_p=(\hat{\mathbb A}^r_{S_p^r},[+])\ra \hat{\mathbb G}_{a,S_p^r}$.
The compositions
$$L_p^r \circ(s_{i,p}-s_{j,p}):\hat{U}_{ij} \hat{\otimes}_{S^0_p} S^r_p \ra \hat{\mathbb G}_{a,S_p^r}$$
define a Cech cocycle of elements
\begin{equation}
\label{vv}
\varphi^r_{ij}\in \cO(\hat{U}_{ij} \hat{\otimes}_{S^0_p} S^r_p)\end{equation}
 and hence a cohomology class $\varphi^r$ in
$H^1(E\hat{\otimes}_{S^0_p} S^r_p,\cO)=H^1(E\otimes_{S^0_p} S^r_p,\cO)$.  The expression
 \begin{equation}
 \label{expression}
 \langle \varphi^r,x \rangle x^{-1-\phi^r},\end{equation}
 where the brackets mean Serre duality, is a well defined element of $S_p^r \cdot x^{-1-\phi^r}$.

 If $L$ is not necessarily free on $X$ we can make the above
 construction locally and the various expressions (\ref{expression}) glue together to give an element
\begin{equation}
\label{cac} f^r_p=f^r_{p,jet}\in M^r_p(-1-\phi^r).\end{equation}

 \subsection{$\d_{\pi}$-overconvergence of $f^r_p$}

 \begin{theorem}
 \label{overfr}
 Assume $v_p(\pi)\geq \frac{1}{p-1}$.
 Then the forms $$\frac{p}{\pi} f^r_p \in M^r_p(-1-\phi^r) \otimes_{R_p} R_{\pi}$$ belong to the image of the homomorphism
 $$M^r_{\pi}(-1-\phi^r)\ra M^r_p(-1-\phi^r) \otimes_{R_p} R_{\pi}.$$
In particular $f^r_p$ are $\d_{\pi}$-overconvergent.
 \end{theorem}

{\it Proof}.
The question is clearly local on $X$ in the Zariski topology so we may assume that $L$ is free on $X$ with basis $x$
and $X$ has an \'{e}tale coordinate $t:X \ra {\mathbb A}^1$. We may also assume that each $U_i \ra X$ factors though
an \'{e}tale map $t_i:U_i\ra X \times {\mathbb A}^1$. Next we note that
\begin{equation}
\label{messs}
\hat{U}_i \hat{\otimes}_{S_p^0} S_p^r \otimes_{R_p} R_{\pi} \simeq (\hat{U}_i \otimes_{R_p} R_{\pi}) \otimes_{S_{\pi}^0} (S_p^r \otimes_{R_p} R_{\pi}).\end{equation}
(This follows from the general fact that if $S$ is a ring, $S', C$ are $S$-algebras, $A,B$ are $C$-algebras,  and $A'=A\otimes_S S'$, $B'=B \otimes_S S'$, $C'=C \otimes_S S'$, then $A \otimes_C B \otimes_S S' \simeq A \otimes_C B'
\simeq A \otimes_C C' \otimes_{C'} B' \simeq A'\otimes_{C'} B'$.) Consequently there is a canonical homomorphism
 from (\ref{messs}) to $\hat{U}_{i,R_{\pi}} \hat{\otimes}_{S_{\pi}^0}S_{\pi}^r$, where, as usual, $\hat{U}_{i,R_{\pi}}=\hat{U}_i \otimes_{R_p} R_{\pi}$.
We claim that one can find sections $s_{i,p}$ and $s_{i,\pi}$ of the canonical projections making the following diagram commute:
\begin{equation}
\label{diaa}
\begin{array}{rcl}
\hat{U}_i \hat{\otimes}_{S_p^0} S_p^r \otimes_{R_p} R_{\pi} & \stackrel{s_{i,p}}{\longrightarrow} & J_p^r(U_i) \otimes_{R_p} R_{\pi}\\
\  & \  & \  \\
 \downarrow & \  & \downarrow  \\
\  & \  & \  \\
 \hat{U}_{i,R_{\pi}} \hat{\otimes}_{S_{\pi}^0}S_{\pi}^r & \stackrel{s_{i,\pi}}{\longrightarrow} & J^r_{\pi}(U_{i,R_{\pi}})
\end{array}\end{equation}
where the vertical morphisms are the canonical ones.
Indeed consider the ring $B=\cO(\hat{U}_{i,R_{\pi}})$ and the commutative diagram
\begin{equation}
\label{diagg}
\begin{array}{rcl}
B[\d_p t,...,\d_p^r t]\h & \leftarrow & B[\d_p t,...,\d_p^r t,\d_p t_i,...,\d_p^r t_i]\h\\
\  & \  & \  \\
\uparrow & \  & \uparrow\\
\  & \  & \  \\
B[\d_{\pi} t,...,\d_{\pi}^r t]\h & \leftarrow & B[\d_{\pi} t,...,\d_{\pi}^r t,\d_{\pi} t_i,...,\d_{\pi}^r t_i]\h \end{array}
\end{equation}
with horizontal arrows sending $\d_p t_i,...,\d_p^r t_i$ and $\d_{\pi} t_i,...,\d_{\pi}^r t_i$ into $0$. Then the spaces in the diagram
(\ref{diaa}) are the formal spectra of the rings in the diagram (\ref{diagg}) and we can take the horizontal arrows in the diagram
(\ref{diaa}) to be induced by the horizontal arrows in the diagram (\ref{diagg}). The diagram (\ref{diaa}) plus Proposition
\ref{xxx}
then induces a commutative
diagram
\begin{equation}
\label{striga}
\begin{array}{rclll}
\hat{U}_{ij} \hat{\otimes}_{S_p^0} S_p^r \otimes_{R_p} R_{\pi} & \stackrel{s_{i,p}-s_{j,p}}{\longrightarrow} &
N_p^r\otimes_{R_p} R_{\pi} & \stackrel{\frac{p}{\pi}L_p^r}{\longrightarrow} & \hat{\mathbb G}_{a,S_p^r\otimes_{R_p} R_{\pi}}\\
\  & \  & \  & \  & \  \\
 \downarrow & \  &  \downarrow  & \  & \downarrow \\
\  & \  & \  & \  & \ \\
\hat{U}_{ij,R_{\pi}} \hat{\otimes}_{S^0_{\pi}} S_{\pi}^r & \stackrel{s_{i,\pi}-s_{j,\pi}}{\longrightarrow} &
N_{\pi}^r & \stackrel{L_{\pi}^r}{\longrightarrow} & \hat{\mathbb G}_{a,S_{\pi}^r}
\end{array}\end{equation}
where $N_{\pi}^r$ is the kernel of the canonical projection $J^r_{\pi}(E_{R_{\pi}})\ra \hat{E}_{R_{\pi}} \hat{\otimes}_{S_{\pi}^0}S_{\pi}^r$ and the vertical morphisms are the canonical ones.
The diagram (\ref{striga}) shows that the cocycle $\frac{p}{\pi}\varphi^r_{ij}$ in (\ref{vv}) comes from a cocycle of elements
in $\cO(\hat{U}_{ij,R_{\pi}} \hat{\otimes}_{S_{\pi}^0}S_{\pi}^r)$. This immediately implies that
the element $\frac{p}{\pi} \langle \varphi^r, x \rangle \in S_p^r \otimes_{R_p} R_{\pi}$ comes from an element in $S_{\pi}^r$ and we are done.
\qed

\begin{remark}
Since $f^1_p\in M^1_p$ is $\d_{\pi}$-overconvergent it follows that its $\d_p$-Fourier expansion $E(f^1_p) \in R_p((q))[\d_p q]\h$ is also $\d_{\pi}$-overconvergent. But, as shown in \cite{difmod},   $E(f^1_p)$ equals the series $\Psi_p$ in (\ref{eqqq}) and note that we knew already (cf. the remarks surounding Equation \ref{eqqq}) that $\Psi_p$ is $\d_{\pi}$-overconvergent. A similar remark holds for $f^r_p$, $r \geq 2$.
\end{remark}

\subsection{$\d_{\pi}$-overconvergence of  $f_p^{\partial},f_{\partial,p}$}
In this subsection we assume that $X \subset X_1(N)_{R_p}$ is an affine open set disjoint from $(cusps)$ and $(ss)$.
There are remarkable forms $f^{\partial}_p \in M^1_p(\phi-1)$
$f_{\partial,p}\in M^1(1-\phi)$
playing a key role in the theory; they are multiplicative inverses of
 each other: $f^{\partial}_p \cdot f_{\partial,p}=1$.
These forms were  introduced in \cite{Barcau} in the level $1$ case; cf. \cite{book}, p. 269, for the arbitrary level case. The definition of these forms in \cite{book}, loc.cit. is crystalline but an alternative description of these forms (up to a multiplicative factor in $R^{\times}$) can be given via \cite{book}, Propositions 8.64 and 8.66; here we shall follow this latter approach. Indeed one has a canonical $R$-derivation
$\partial:\cO(V)\ra \cO(V)$ defined by Katz \cite{Katz} via the Gauss-Manin connection, generalizing the ``Serre operator"; cf. \cite{book}, pp.254-255, for a review of this. (Here $V$ is as in (\ref{defV}).)
One can consider then the conjugate operators $\partial_0, \partial_1:M^1_p \ra M^1_p=\cO(J^1_p(V))$; cf. (\ref{cirip}).
One can also consider the {\it Ramanujan form} $P \in M_p^0(2)$; cf. \cite{book}, p. 255, for a review of this. Then one can define $f^{\partial}_p \in M^1_p$ by the formula
\begin{equation}
\label{deffpartial}
f^{\partial}_p:=\partial_1 f^1_p-p\phi(P)f^1_p \in M^1_p.
\end{equation}
It turns out that actually $f^{\partial}_p$ has weight $\phi-1$, i.e. $f^{\partial}_p\in M^1_p(\phi-1)$ and $f^{\partial}_p$ has a multiplicative inverse in $M^1_p$ called $f_{\partial,p}$ which actually lies in $M^1(1-\phi)$ and satisfies
\begin{equation}
\label{foo}
f_{\partial,p}=-\partial_0 f^1_p+P \cdot f^1_p.
\end{equation}
 (By the way, as shown in \cite{Barcau}, $f^{\partial}_p$ and hence $f_{\partial,p}$, have $\d_p$-Fourier expansion $E(f^{\partial}_p)=E(f_{\partial,p})=1$.)

Theorem \ref{overfr} plus Proposition \ref{conjugate} imply then the following:

\begin{theorem}
\label{overfpartial}
Assume $v_p(\pi)\geq \frac{1}{p-1}$.
Then the elements $\frac{p}{\pi} f^{\partial}_p, \frac{p}{\pi} f_{\partial,p} \in M^1_p\otimes_{R_p} R_{\pi}$ belong to the image of the map $M^1_{\pi}\ra M^1_p\otimes_{R_p} R_{\pi}$. In particular $f^{\partial}_p$ and $f_{\partial,p}$ are $\d_{\pi}$-overconvergent.
\end{theorem}

\subsection{Review of the $\d_p$-characters $\psi_p$ of elliptic curves \cite{book}}
We follow \cite{book}, pp. 194-197. Let $A/R_p$ be an elliptic curve  and fix  a level $\Gamma_1(N)$ structure on $A$.
(The construction below does not depend on this level structure.)
If  $Y_1(N)_{R_p}:=X_1(N)_{R_p}\backslash (cusps)$ we get an induced point $P_A:Spec\ R_p\ra Y_1(N)_{R_p}$. Let $X \subset Y_1(N)_{R_p}$ be an affine open set ``containing" the above point and such that the line bundle $L$ on $X$ is trivial with basis $x$. Let $\omega$ be the invertible $1$-form on $A$
defined by $x$. By the universality property of the $p$-jet spaces we get canonical morphisms $P^r_A:\cO(J^r_p(X))\ra R$ compatible with $\d_p$ in the obvious sense. Then any $\d_p$-modular form $f \in M^r_p$ on $X$ defines an element $f(A,\omega)\in R_p$ as follows: we write $f=\tilde{f} \cdot x^w$ with $\tilde{f} \in \cO(J^r_p(X))$ and one takes $f(A,\omega)\in R_p$ to be the image of $\tilde{f}$ in $R_p$
via the above morphism $P^r_A$. In particular one can consider the $\d_p$-modular forms
$f^1_p\in M^1_p(-1-\phi)$ and $f^2_p\in M^2_p(-1-\phi^2)$ in (\ref{cac}); we get elements
$f^1_p(A,\omega), f^2_p(A,\omega)\in R_p$. We recall that $f^1_p(A,\omega)=0$ if and only if $A$
has a lift of Frobenius i.e. the $p$-power Frobenius of $A\otimes_{R_p} k$ lifts to a morphism of schemes $A \ra A$ over ${\mathbb Z}$. Assume in what follows that {\it $A$ does not have lift of Frobenius}.
Then the quotient $\frac{f^2_p(A,\omega)}{f^1_p(A,\omega)}$, which is a priori an element of $K_p$, lies actually in $R_p$.
On the other hand we may consider the cocycles (\ref{vv}).
The images of these cocycles via the homomorphism $S_p^r=\cO(J^r_p(X))\ra \cO(J^r_p(X'))\stackrel{P^r_A}{\ra}R_p$ yield cocycles
$$\varphi^r_{ij}(A)\in \cO(\hat{U}_{ij,A})$$
where $U_{ij,A}=U_{ij}\cap A$. (Here we view $A$ embedded into the universal elliptic curve $E$ via the isomorphism $A\simeq E \times_{X,P_A} R_p$.) The cocycle
$$\varphi^2_{ij}(A)-\frac{f^2_p(A,\omega)}{f^1_p(A,\omega)} \varphi^1_{ij}(A) \in \cO(\hat{U}_{ij,A})$$
turns out, by construction, to be a coboundary $$\Gamma_i-\Gamma_j$$
with $\Gamma_i\in \cO(\hat{U}_{i,A})$, $U_{i,A}=U_i \cap A$.
Recall the series $L^r_p \in S_p^r[\d_p T,..., \d_p^r T]\h$; cf. (\ref{Lpr}). (Here $T$ is an \'{e}tale coordinate at the origin of $E$ such that $x\equiv dT$ mod $T$.) The images of $L^r_p$  via $S_p^r \ra R_p$ yield series $L^r_p(A) \in R_p[\d_p T,...,\d_p^r T]\h$. Take sections $s_{i,p}:\hat{U}_{i,A} \ra J^2_p(U_{i,A})$
of the natural projections and let $N^2_{p,A}$ be the kernel of the projection $J^2_p(A)\ra \hat{A}$.
The maps
\begin{equation}
\label{taui}
\tau_{i,p}:\hat{U}_{i,A} \hat{\times} N^2_{p,A} \ra J^2_p(U_{i,A}),\end{equation}
 given at the level of points
by $(a,b)\mapsto s_{i,p}(a)+b$, are isomorphisms. Consider the functions
\begin{equation}
\label{hopper}\begin{array}{rcl}
\psi_{i,p} & := & L^2_p(A)-\frac{f^2_p(A,\omega)}{f^1_p(A,\omega)} L^1_p(A)+\Gamma_i\\
\  & \  & \  \\
 \  & \in & \cO(\hat{U}_{i,A})[\d_p T,\d_p^2 T]\h\\
 \ & \  & \  \\
 \  &  = & \cO(\hat{U}_{i,A}\hat{\times} N^2_{p,A}). \end{array}\end{equation}
Then it turns out that the functions
$$\psi_{i,p}\circ \tau_{i,p}^{-1}\in \cO(J^2_p(U_{i,A}))$$
glue together to give a function
\begin{equation}
\label{deltacharacter}
\psi_p\in \cO(J^2(A)).\end{equation}
This map turns out to be an homomorphism $J^2_p(A)\ra \hat{\mathbb G}_a$ and was referred to in \cite{book}, Definition 7.24, as the {\it canonical} $\d_p$-character (of order $2$) of $A$. (In loc. cit. $\psi_p$ was denoted by $\psi_{can}$.)

In case $A$ has a lift of Frobenius a different (but similar, and in fact easier) construction leads to what in cf. \cite{book}, Definition 7.24 was referred to as the  {\it canonical
$\d_p$-character} (of order $1$)  of $A$.
We will denote it again by
\begin{equation}
\label{woo}\psi_p \in \cO(J^1_p(A)).\end{equation} In \cite{book}, loc.cit. this $\d_p$-character was again denoted by $\psi_{can}$.

\subsection{$\d_{\pi}$-overconvergence of $\psi_p$}

Let $A/R_p$ be an elliptic curve and let $r$ be $1$ or $2$ according as $A$ has a lift of Frobenius or not.

\begin{theorem}
\label{overpsi}
Assume $v_p(\pi)\geq \frac{1}{p-1}$. Then
the function $\frac{p}{\pi}\psi_p$ belongs to the image of the map
$$\cO(J^r_{\pi}(A_{R_{\pi}})) \ra \cO(J^r_p(A))\otimes_{R_p}R_{\pi}.$$
In particular $\psi_p$ is $\d_{\pi}$-overconvergent.
\end{theorem}

{\it Proof}. We give the proof in case $r=2$. The proof in case $r=1$ is similar. It is enough to show that one can choose the data in our construction such that:

1) The functions $\frac{p}{\pi} \psi_{i,p}$ (where $\psi_{i,p}$ is as in (\ref{hopper})) belong to the image of
$$\cO(\hat{U}_{i,A}\otimes_{R_p}R_{\pi})[\d_{\pi} T, \d_{\pi}^2T]\h \ra \cO(\hat{U}_{i,A}\otimes_{R_p}R_{\pi})[\d_p T, \d_p^2T]\h;$$

2) There are commutative diagrams
$$\begin{array}{rcl}
(\hat{U}_{i,A}\otimes_{R_p}R_{\pi})\hat{\times} (N^2_{p,A}\otimes_{R_p}R_{\pi}) & \stackrel{\tau_{i,p}}{\ra} &  J^2_p(U_{i,A})\otimes_{R_p}R_{\pi}\\
\downarrow & \  & \downarrow\\
(\hat{U}_{i,A}\otimes_{R_p}R_{\pi}) \hat{\times} N^2_{\pi,A} & \stackrel{\tau_{i,\pi}}{\ra} &
J^2_{\pi}(U_{i,A}\otimes_{R_p}R_{\pi})\end{array}$$
for isomorphisms $\tau_{i,\pi}$.

Now 1) follows from the fact that  $\frac{p}{\pi}L_p^r(A) \in R_{\pi}[\d_{\pi} T,\d_{\pi}^2 T]\h$ (cf. Theorem \ref{xxx}), and $\Gamma_i \in \cO(\hat{U}_{i,A})$.
On the other hand 2) follows from the fact that one can choose the sections $s_{i,p}$ together with sections $s_{i,\pi}$ as in (\ref{diaa}); then one can define the isomorphisms $\tau_{i,\pi}$ using
$s_{i,\pi}$ in the obvious way.
This ends the proof.
\qed

\subsection{$\d_{\pi}$-overconvergence of $f^{\sharp}_p$ for $f$ on $\Gamma_0(N)$}\
We first recall the construction of the $\d_p$-modular forms $f^{\sharp}_p$ attached to newforms on $\Gamma_0(N)$ given in \cite{eigen,BP}. As usual we let $N>4$, $(N,p)=1$.
 Fix, in what follows, a normalized newform
$f=\sum_{n \geq 1} a_n q^n$
 of weight $2$ on $\Gamma_0(N)$ over ${\mathbb Q}$ and an  elliptic curve $A$ over ${\mathbb Q}$ of conductor $N$ such that $f$ and $A$ correspond to each other in the sense of Theorem \ref{fundamental}; recall that this means that
 there exists
 a morphism
\begin{equation}
\label{shimuraa}
\Phi:X_0(N)\ra A
 \end{equation}
 over ${\mathbb Q}$  such that the pull back to $X_0(N)$ of some $1$-form on $A$ over ${\mathbb Q}$
corresponds to $f$ and $L(A,s)=\sum a_n n^{-s}$.
Fix  an embedding ${\mathbb Z}[1/N,\zeta_N]\subset R_p$. Let $A_{R_p}$ be the N\'{e}ron model of $A\otimes_{\mathbb Q} K_p$ over $R_p$ (which is an elliptic curve) and let $X_1(N)_{R_p}$ be the (smooth)  ``canonical" model of $X_1(N)$ over $R_p$ which has been considered before. By the N\'{e}ron model property there is an induced morphism $\Phi_p:X_1(N)_{R_p}\ra A_{R_p}$. Let $X \subset X_1(N)_{R_p}$ be any affine open set.
Let $r$ be $1$ or $2$ according as $A_{R_p}$ has or has not a lift of Frobenius. (Note that we always have $r=2$ if $A$ has no complex multiplication.)
The image of the canonical $\d_p$-character $\psi_p
\in \cO(J^r_p(A_{R_p}))$
in (\ref{deltacharacter}) (respectively (\ref{woo})) via the map
$$\cO(J^r_p(A_{R_p})) \stackrel{\Phi_p^*}{\longrightarrow} \cO(J^r_p(X))=S^r_p=M^r_p(0)\subset M^r_p$$
is denoted by $f^{\sharp}=f^{\sharp}_p$ and is a $\d_p$-modular form of weight $0$; this form was introduced in \cite{eigen} and played a key role in \cite{BP}.

Putting together Theorem \ref{overpsi} and Remark \ref{usefull} we get:

\begin{theorem}
\label{lastover} Assume $v_p(\pi)\geq \frac{1}{p-1}$. Then
the function $\frac{p}{\pi}f^{\sharp}_p$ belongs to the image of the map
$$\cO(J^r_{\pi}(X_{R_{\pi}})) \ra \cO(J^r_p(X))\otimes_{R_p}R_{\pi}.$$
In particular $f^{\sharp}_p$ is $\d_{\pi}$-overconvergent.
\end{theorem}

\bibliographystyle{amsplain}

\end{document}